\documentclass[12pt]{article}
\usepackage[cp1251]{inputenc}
\usepackage[russian]{babel}
\begin{document}
\author{Белошапка В.К.}

\date{}

\title{Модификация конструкции Пуанкаре\\ и её применение в $CR$-геометрии\\ гиперповерхностей в ${\bf C}^4$}

\maketitle

\begin{abstract}
Обобщение гомологического оператора Пуанкаре -- модифицированная конструкция Пуанкаре -- была использована для оценки размерности локальных автоморфизмов произвольного ростка вещественно аналитической гиперповерхности пространства $\mathbf{C}^4$. В работе доказана следующая альтернатива: либо эта размерность  бесконечна, либо она не превосходит 24-х. При этом 24 реализуется лишь для невырожденной гиперквадрики (одной из двух). Если гиперповерхность 2-невырождена в точке общего положения,  то оценку можно улучшить до 17, а если  3-невырождена, то до 20.
\end{abstract}

\footnote{
Механико-математический факультет Московского университета им.Ломоносова,
Воробьевы горы, 119992 Москва, Россия, vkb@strogino.ru }

{\bf Введение}

\vspace{3ex}

Ведущим элементом метода модельной поверхности  является конструкция Пуанкаре, которую он применял как
в небесной механике, так и в $CR$-геометрии (гомологический оператор Пуанкаре или Пуанкаре-Дюлака).   Применение в $CR$-геометрии дано в работе 1907-го года \cite{P}. Эта конструкция, по существу, представляет собой версию теоремы о неявном отображении в классе формальных степенных рядов.

Обычно в $CR$-геометрии эта конструкция используется так. Пусть имеется некторое нелинейное дифференциальное или функциональное соотношение $F(x,\phi(x))=0$.  Пусть в кольце формальных степенных рядов от $x$ введена некоторая градуировка (вес), причем $\mu$-я компонента нашего соотношения имеет вид
$$L(x,\phi_{\mu}(x)) =  \mbox{ выражению, зависящему от   }   \phi_{\nu}  \mbox{  при   }  \nu  \leq \mu-1,$$
где  $L(x,y)$  линейно зависит от $y$. Тогда, очевидно, размерность линейного пространства решений уравнения $L(x,\phi(x))=0$
 мажорирует размерность семейства решений исходного уравнения  $F(x,\phi(x))=0$.

      В работе \cite{VB05} для оценки размерности алгебры Ли инфинитезимальных автоморфизмов произвольной голоморфно невырожденной вещественной гиперповерхности пространства $\mathbf{C}^3$ была использована некоторая модификация этой конструкции (редукция на глубину два). А именно, пусть  $\mu$-я компонента нашего соотношения $F(x,\phi(x))=0$ имеет вид
$$L_1(x,\phi_{\mu}(x))+L_2(x, \phi_{\mu -1}(x)) =  \mbox{  выражению, зависящему от   }   \phi_{\nu}  \mbox{  при   }  \nu \leq \mu-2,$$
где  $L_1(x,y)$ и $L_2(x,y)$  линейно зависят от $y$. Тогда, очевидно, размерность линейного пространства решений уравнения $L_1(x,\phi(x))+L_2(x,\phi(x))=0$  мажорирует размерность семейства решений исходного уравнения  $F(x,\phi(x))=0$. Аналогично определяется обобщение этой конструкции -- редукция на произвольную глубину $k$.

      В данной работе мы даем еще одну демонстрацию применения редукции на глубину два и три.  С помощью этой модификации конструкции Пуанкаре дается оценка сверху на размерность   алгебры Ли инфинитезимальных автоморфизмов произвольной голоморфно невырожденной вещественной гиперповерхности пространства $\mathbf{C}^4$ (теорема 16).

      Этот результат подтверждает старую гипотезу \cite{VB96}: либо размерность группы  автоморфизмов произвольного ростка вещественно аналитической гиперповерхности не превышает размерности для любой невырожденной  стандартной гиперквадрики (в $\mathbf{C}^4$ она равна 24), либо она бесконечна.

Отметим, что для получения известной оценки для гиперповерхностей пространства $\mathbf{C}^2$ достаточно обычной (однократной) конструкции Пуанкаре.  Для получения же такой оценки в пространстве размерности $(n+1)$ потребуется использование всех вариантов редукции на глубину от 1 до $n$.

\vspace{5ex}

Поскольку $k$-кратная конструкция Пуанкаре отличается от классической,  изложим схему её  применения в общем виде.  Пусть
$V$ - линейное пространство бесконечных последовательностей вещественных чисел,  пусть $x \in V$. Пусть, далее,
эта последовательность разбита на конечные отрезки, которые мы будем обозначать $x_j$. Соответственно $x=(x_1,x_2, \dots)$,
при этом $x_j$ -- элемент  некоторого конечномерного вещественного линейного пространства. Фиксируем некоторое натуральное
число  $k$.  Вот общая формулировка,  которую естественно назвать схемой рекурсии на глубину $k$.

\vspace{2ex}

{\bf Теорема 1:}  Пусть имеется бесконечная система полиномиальных соотношений вида
\begin{eqnarray} \label{T1}
\Theta_j (x_1, \dots, x_j) = L_{j1}(x_j)+ \dots +L_{jk}(x_{j-k+1})+\theta_j(x_{j-k},x_{j-k-1},\dots,x_1)=0,\\
\nonumber
 j=k,(k+1), \dots,    \quad  \mbox{  где  }   (L_{j1},\dots,L_{jk})   \mbox{  -- линейны,   причем  } \\
\nonumber L_j(x)=L_{j1}(x_j)+L_{j2}(x_{j-1})+ \dots +L_{jk}(x_{j-k+1}) \mbox{ и } L(x)=(L_1(x), L_2(x), \dots ).
\end{eqnarray}
И пусть известно, что  ${\rm Ker }\, L$    содержится в конечномерном подпространстве  пространства $V$ вида $\tilde{V}_l=\{(x_1, \dots,x_l,0,0,\dots)\}$. Тогда число параметров, от которых зависит общее решение (\ref{T1}),   не превосходит размерности ${\rm Ker} \, L$.\\
{\it Доказательство:} Пусть $W$ -- прямое дополнение ${\rm Ker} \, L$ до $V$. Это дополнение можно определить так.
В конечномерном пространстве  $\tilde{V}_l$  произвольно выберем прямое дополнение $\tilde{W}$ к $ {\rm Ker} \, L$ и дополним его подпространством $V_l=\{(0,\dots,0,x_{l+1},x_{l+2},\dots )\} $. Таким образом, уравнение $L(x)=0$ имеет в пространстве $W$ единственное решение $x=0$.
Просматривая последовательно уравнения $\Theta_j(x)=L_j(x)+\theta_j(x)=0$, убеждаемся, что эта система также имеет в $W$ не более одного решения. Произвольный вектор $x \in V$ имеет вид $x+a$, где $x \in W, \; a \in   {\rm Ker} \, L$.  Рассмотрим нашу систему при фиксированном $a$.   Получаем $\Theta_j(x+a)=L_j(x)+\theta_j(x+a)=0$.   Также видим, что эта система имеет не более одного решения. Таким образом, совокупность решений (\ref{T1}) параметризуется некторым подмножеством ${\rm Ker} \, L$. Теорема доказана.

\vspace{5ex}

Пусть $\Gamma$ -- вещественно аналитическая гиперповерхность в области пространства $\mathbf{C}^4$,   $\Gamma_{\xi}$ -- росток этой гиперповерхности в
точке $\xi$. Пусть, далее,  ${\rm aut} \, \Gamma_{\xi}$ -- алгебра Ли, состоящая из ростков вещественных вещественно
аналитических векторных полей в точке $\xi$, касательных к  $\Gamma_{\xi}$.    Если $\Gamma$ вне собственного аналитического
подмножества Леви-невырождена, то в рамках стандартного подхода (однократная редукция) мы получаем стандартную оценку:
${\rm dim \; aut} \,\Gamma_{\xi}$   не превосходит размерности алгебры автоморфизмов касательной невырожденной гиперквадрики, которая, независимо от сигнатуры, равна 24.

Рассчитывать на получение оценки мы можем только в случае, если $\Gamma$ -- голоморфно невырождена.  Для гиперповерхности в $\mathbf{C}^4$ голоморфная невырожденность эквивалентна $l$-невырожденности вне собственного аналитического подмножества,
причем $l \leq 3$ (см.\cite{BER}).   Таким образом, для получения общей оценки нам необходимо рассмотреть две разных ситуации. Первая:  $\Gamma$  равномерно 2-невырождена в окрестности ${\xi}$.   Вторая: $\Gamma$ равномерно 3-невырождена в окрестности ${\xi}$.

Примеров равномерно 2-невырожденных гиперповерхностей в $\mathbf{C}^4$ достаточно много.  Например, они содержатся в известной работе Г.Фелса и В.Каупа \cite{FK08}.  В частности, там описаны трубчатые гиперповерхности над вещественными конусами и их голоморфные автоморфизмы. В  $\mathbf{C}^4$ группы таких конусов имеют размерность 15.   Примеров равномерно 3-невырожденных  гиперповерхностей в $\mathbf{C}^4$ нам известно два. Один от В.Каупа \cite{FK08} и один от А.Санти \cite{AS17}.  Оба -- голоморфно однородные, причем в примере Санти известна размерность автоморфизмов, которая равна 8.

\vspace{3ex}

Рассмотрим оба случая -- 3-невырожденный и  2-невырожденный -- последовательно. Отметим, что для получения оценки нам пришлось не только использовать гомологический оператор на глубину больше чем один.   Также для анализа 2-невырожденных  мы использовали двухкратную процедуру.  А именно, рекуррентная процедура с одним весом, затем -- смена веса и анализ ядра старого гомологического оператора с точки зрения новой весовой рекурсии, что ведет к новому гомологическому оператору. Затем для анализа специального случая (большое ядро) -- еще одна смена веса и новая рекурсия.

     В итоге мы получим две оценки размерности. Для 2-невырожденных - 17, а для 3-невырожденных - 20. При этом
для 2-невырожденных  наша техника дает оценку 18, но использование недавнего результата И.Зеленко и Д.Сайкса
\cite{ZS21} позволило улучшить ее до 17.

\vspace{10ex}

{\bf 3-невырожденные гиперповерхности}

\vspace{3ex}

Обозначим координаты в $\mathbf{C}^4$ через $(z, \zeta, \eta, w=u+i \, v)$.
Пусть  $\Gamma$ -- {\it равномерно 3-невырождена} в окрестности ${\xi} \in  \Gamma$.
Имея ввиду нашу цель -- получение оценки на размерность группы автоморфизмов, -- мы можем ограничиться так называемыми "жесткими" $\,$ гиперповерхностями. Действительно, если алгебра автоморфизмов содержит поле, трансверсальное комплексной касательной, то после локального распрямления мы можем полагать, что группа содержит сдвиги вдоль оси $u$ и, тем самым, локальное уравнение  $ \Gamma$
имеет вид
$$ v =F(z,\bar{z},\zeta,\bar{\zeta},\eta,\bar{\eta}),$$
т.е. правая часть не зависит от $u$.
В силу 3-невырожденности ранг формы Леви в общей точке равен 1, и мы можем записать локальное уравнение $\Gamma$ в виде $ v =|z|^2 +F_3+F_4+ \dots,$ где $F_j(z,\bar{z},\zeta,\bar{\zeta},\eta,\bar{\eta})$ -- это однородный полином степени $j$. При этом простыми преобразованиями мы можем удалить все плюригармонические компоненты правой части уравнения, а также все члены, линейные по $z$ и $\bar{z}$, за исключением $|z|^2$.

Необходимым условием равномерной 3-невырожденности является условие, что ранг комплексного гессиана $F(z,\bar{z},\zeta,\bar{\zeta},\eta,\bar{\eta})$ всюду не превосходит единицы. Это условие, учитывая, что $F_{z \bar{z}}$ в нуле равна единице , можно записать как условие равенства нулю трех миноров второго порядка. А именно
\begin{eqnarray} \label{3ND}
\nonumber
 \delta_1(F)=F_{z \bar{z}} \,  F_{\zeta \bar{\zeta}}- | F_{z \bar{\zeta}}|^2=0,\\
\nonumber  \delta_2(F)=F_{z \bar{z}} \,  F_{\zeta \bar{\eta}}-F_{\zeta \bar{z}} \,  F_{z \bar{\eta}}=0,\\
 \delta_3(F)=F_{z \bar{z}} \,  F_{\eta \bar{\eta}}- | F_{z \bar{\eta}}|^2=0.
\end{eqnarray}

\vspace{3ex}

{\bf Лемма 2:}    Если $\Gamma$  равномерно 3-невырождена, то уравнение этой гиперповерхности после линейной замены можно записать в виде
\begin{eqnarray} \label{J6}
v=|z|^2+ F_3+F_4 +F_5+F_6+O(7),   \qquad \qquad  \qquad \qquad\\
\nonumber
 \mbox{где} \quad F_3=2 \, {\rm Re} (z^2 \, \bar{\zeta}), \;  F_4=2 \, {\rm Re} (z^3 \, \bar{\eta})+4 \,|z|^2 \, |\zeta|^2,  \qquad  \qquad\\
\nonumber
F_5=2 \, {\rm Re} (r_{{1}}{{\bar{z}}}^{4}\zeta+ r_{{2}}{{\bar{z}}}^{4}\eta+ r_{{3}}z{{\bar{\zeta}}}^{2}{{\bar{\eta}}}^{2}+r_{{4}}z{{\bar{\eta}}}^{4}+
4\,{z}^{2}\zeta{{\bar{\zeta}}}^{2}+6\,{z}^{2}\bar{z}\,\zeta{\bar{\eta}}), \qquad\\
\nonumber
F_6=2 \, {\rm Re} (
8\,\bar{r}_{{1}}{z}^{3}{\bar{z}}\, \zeta \,{\bar{\zeta}}+ 8\, \bar{r}_{{2}}{z}^{3}{\bar{z}}\,\zeta \,{\bar{\eta}}+2\,r_{{3}}z{\eta}^{2}{\bar{\zeta}}
\,{\zeta}^{2}+2\,r_{{4}}z{\eta}^{4}{\bar{\zeta}}+ \qquad \\
\nonumber s_{{1}}z{{\bar{z}}}^{4}\, \zeta+s_{{2}}{{\bar{z}}}^{5} \, \zeta +s_{{3}}z{{\bar{z}}}^{4} \, \eta+s_{{4}}{\bar{z}}^{5} \eta+s_{{5}}{{\bar{z}}}^{4}{\zeta}^{2}+ s_{{6}}{{\bar{z}}}^{4} \, \zeta \, \eta +s_{{7}}{{\bar{z}}}^{4}{\eta}^{2}+
 s_{{8}}\bar{z}{{\eta}}^{5}+\\
\nonumber  12\,{{\bar{z}}}^{3} \zeta{\bar{\zeta}}\, \eta+12\,z{{\bar{z}}}^{2}{\zeta}^{2}{\bar{\eta}}) + 16\,|z|^2\,|{\zeta}|^{4}+9\,|{z}|^{4} \, |\eta|^2  \qquad \qquad \qquad
\end{eqnarray}
{\it Доказательство:}  Запишем общий вид $F_3$ с учетом наших упрощений. Выделяя в соотношениях (\ref{3ND}) компоненту степени один, получаем, что  $F_3=2 \, {\rm Re} (a_1 \, \zeta + a_2 \, \eta) \, \bar{z}^2 $. Поскольку $F_3$  не есть тождественный ноль, то $(a_1,a_2) \neq 0$ и мы можем заменить $ {a}_1 \, \zeta +{a}_2 \, \eta$ на новую переменную $\zeta$.  Теперь простым преобразованием мы можем удалить из всех последующих компонент $F$ слагаемые вида $ 2 \, {\rm Re} (A(z,\zeta,\eta) \, \bar{z}^2)$  с голоморфным коэффициентом $A$.

 Запишем общий вид $F_4$ с учетом наших упрощений. Выделяя в соотношениях  (\ref{3ND})  компоненту степени два, получаем, что
$F_4=2 \, {\rm Re}\, \bar{z}^3 \, (a_3 \, \eta+ a_4 \zeta)+4 \,|z|^2 \, |\zeta|^2$. Из равномерной 3-невырожденности следует, что   $(a_3,a_4) \neq 0$ и мы можем заменить $ {a}_3 \, \eta +{a}_4 \, \zeta$ на новую переменную $\eta$. Простым преобразованием удаляем из всех последующих компонент $F$ слагаемые вида $ 2 \, {\rm Re} (B(z,\zeta,\eta) \, \bar{z}^3)$ с голоморфным коэффициентом $B$.
Выделяя компоненту степени три, получаем указанный вид $F_5$. А затем из компоненты степени четыре -- вид $F_6$.
 Лемма доказана.

\vspace{3ex}

Рассмотрим отображение ростка одной гиперповерхности $\Gamma$  в начале координат вида  (\ref{J6}) на другую такую гиперповерхность $\tilde{\Gamma}$. Пусть  координаты ростка  отображения в начале координат имеют вид
$$\Phi=(z  \rightarrow  f(z,\zeta,\eta,w),  \; \zeta  \rightarrow g(z,\zeta,\eta,w),  \; \eta  \rightarrow h(z,\zeta,\eta,w),  \;w  \rightarrow e(z,\zeta,\eta,w)).$$
Будем считать эти гиперповерхности фиксированными.
Введем в пространстве степенных рядов от $(z,\bar{z},\zeta, \bar{\zeta},\eta, \bar{\eta},u)$, а также от
$(z,\bar{z},\zeta,\bar{\zeta}, \eta, \bar{\eta}, w,\bar{w})$ градуировку, назначая веса переменным
$$ [z]=[\bar{z}]=[\zeta]=[ \bar{\zeta}]=[\eta]=[ \bar{\eta}]=1, \; [w]=[\bar{w}]=[u]=2.$$
Набор весовых компонент $(f_{\mu-1},g_{\mu-2},h_{\mu-3},e_{\mu})$ обозначим через $\phi_{\mu}$ ( $\mu$-я весовая компонента $\Phi$). Запишем соотношение, отражающее тот факт, что $\Phi$ отображает  $\Gamma$ на $\tilde{\Gamma}$.
\begin{eqnarray} \label{mid}
 \nonumber
  \Theta(z,\bar{z},\zeta, \bar{\zeta},u)=- 2 \, {\rm Im} \, e(z,\zeta,\eta,w)+2 \, |f|^2+4 \, {\rm Re}(f^2 \bar{g}) +4 \, {\rm Re}(f^3 \bar{h}) + 8 \,|f|^2 \, |g|^2+\\
 \nonumber  2 \, F_4(f,\bar{f},g, \bar{g},h,\bar{h}) +2 \, F_5(f,\bar{f},g, \bar{g},h,\bar{h}) +2 \, F_6(f,\bar{f},g, \bar{g},h,\bar{h})   +\dots\\
 \mbox{   при }    w=u+i \, (|z|^2 + 2 \, {\rm Re} (z^2 \, \bar{\zeta})+2 \, {\rm Re} (z^3 \, \bar{\eta})+4 \,|z|^2 \, |\zeta|^2 + \dots)
  \end{eqnarray}

Среди всех голоморфных в начале координат отображений  выделим класс отображений вида
\begin{equation}\label{V}
\mathcal{V}_5=\{ \Phi=Id+\phi_5+ \dots=(z+O(4),  \;  \zeta+O(3) ,  \; \eta +O(2),  \;w+O(5)) \}
\end{equation}

Оценку размерности семейства таких отображений $\Gamma$ на  $\tilde{\Gamma}$ проведем по схеме кратной рекурсии глубины $k=3$ (см. теорему 1). Для этого выделим в $\Theta_{\mu}$ -- $\mu$-й весовой компоненте (\ref{mid}) -- члены, зависящие от $(\phi_{\mu}, \phi_{\mu-1},\phi_{\mu-2})$, т.е.
от
$$(e_{\mu}, e_{\mu-1},e_{\mu-2}, f_{\mu-1}, f_{\mu-2},f_{\mu-3},g_{\mu-2}, g_{\mu-3},g_{\mu-4},h_{\mu-3}, h_{\mu-4},h_{\mu-5}).$$
Введем обозначения  $\Delta_1 \, \psi(u) = i  \, F_3  \, \psi'(u), \; \Delta_2 \, \psi(u) = i  \, F_4  \, \psi'(u).$
Выделим последовательно в слагаемых выражения $\Theta$,  начиная от $- 2 \, {\rm Im} \, e$ до $F_6$, члены указанного вида, получим следующий результат.

\vspace{2ex}

{\bf Лемма 3:}   Для всех $\mu \geq 5$  $\mu$-я компонента $\Theta$ имеет вид
\begin{eqnarray*}
\nonumber
\Theta_{\mu}=L_1(\phi_{\mu})+ L_2(\phi_{\mu-1})+L_3(\phi_{\mu-2})+ \theta_{\mu}(\phi_{\nu < \mu-2}),
 \mbox{  где }  w=u+i \,|z|^2, \\
L_1(\phi) = 2 \, {\rm Re}(i \, e + 2 \, \bar{z} \, f + 2 \, \bar{z}^2 \,g + 2 \, \bar{z}^3 \,h),\\
L_2(\phi) =  \Delta_1 (L_1(\phi))  +  l_2(\phi) , \quad   L_3(\phi) =  \Delta_1 (L_2 (\phi))  +\Delta_2 (L_1(\phi))  +  \Delta^2_1 (L_1(\phi))  + l_3(\phi) ,\\
l_2(\phi)= 2 \, {\rm Re} \{
4 \,z \, \bar{\zeta} \,f+8 \,z \, \bar{z} \,  \bar{\zeta} \,\eta \, g + (4 \, \bar{z}^4 \, r_2 + 4 \, \bar{z} \, \zeta^2 \, \eta \, r_3 + 8 \,\bar{z} \,\eta^3 \, r_4 + 12 \,  \bar{\zeta} \, \bar{z}^2 \,z) \, h \},\\
l_3(\phi)= 2 \, {\rm Re}
 \{(8\,\bar{\zeta}\,\bar{z}\,\zeta+6\,\bar{\eta}\,z^2)\,f+(4\,\bar{z}^4\,r_1+4\,\bar{z}\,\zeta^2\,\eta\,r_3+12\,\bar{\eta}\,\bar{z}\,z^2+\\
8\,\bar{\zeta}^2\,z+16\,\bar{\zeta}\,\bar{z}\,\zeta)\,g+
(16\,\bar{\zeta}\,\bar{z}^3\,z\,r_2+8\,\bar{\zeta}\,\zeta^2\,\eta\,z\,r_3+16\,\bar{\zeta}\,\eta^3\,z\,r_4+ 2\,\bar{z}^5\,s_4+2\,\bar{z}^4\,\zeta\,s_6+\\
4\,\bar{z}^4\,\eta\,s_7+2\,\bar{z}^4\,z\,s_3+10\,\bar{z}\,\eta^4\,s_8+24\,\bar{\zeta}^2\,cz^2\,z +
\nonumber 24\,\bar{\zeta}\,\bar{z}^3\,\zeta+18\,\bar{\eta}\,\bar{z}^2\,z^2)\,h \}.
\end{eqnarray*}

\vspace{3ex}

Отметим, что выражение $L (\phi)= L_1 (\phi)+ L_2(\phi)+ L_3(\phi)$  линейно по $\phi$ и не зависит от $\mu$.

\vspace{3ex}

Пусть  $V_5$ - линейное пространство, состоящее из ростков формальных степенных рядов в начале координат вида
$$\Phi = \phi_5  +  \phi_6 +\dots = (f_4 +f_5+\dots, g_3 +g_4+\dots, h_2 +h_3+\dots, e_5+e_6+\dots )$$
Тогда, в соответствии с теоремой 1, размерность семейства отображений    $\Gamma$ на  $\tilde{\Gamma}$ из $\mathcal{V}_5$
не превосходит размерности ядра $L$ на $V_5$.

\vspace{3ex}

Перейдем к оценке размерности ядра оператора $L$, т.е. пространства решений соотношения
\begin{equation}\label{ker}
  L(\phi)=0, \quad  \mbox{где}   \; \phi \in  V_5
\end{equation}

Обозначим
$$f(0,0,0,u)=a(u),  \; g(0,0,0,u)=b(u), \; h(0,0,0,u)=c(u), \; e(0,0,0,u)=d(u).$$
Положим в (\ref{ker}) $ \bar{z}= \bar{\zeta}=\bar{\eta}=0$, получим соотношение, из которого сразу находим
\begin{eqnarray} \label{e}
 e(z,\zeta,\eta,u)=    \qquad \qquad\qquad\qquad\qquad\qquad\qquad\qquad\qquad\qquad\qquad\\
\nonumber  d(u)+ 2 \, i  \, z  \, \bar{a}(u) +   2 \, i  \, z^2  \, \bar{b}(u) + 2  \, i  \, z^3  \, \bar{c}(u) +
 4 \, i  \, z^4  \, (\bar{r}_2  \, \bar{c}(u)+ \bar{r}_1  \bar{b}(u))+
 2  \, i  \, \bar{s}_4  \, z^5  \, \bar{c}(u),
\end{eqnarray}
причем $d(u)$ -- вещественна.

  Обозначим
\begin{eqnarray*}
f'_z(0,0,0,u)=k_1(u),  \; g'_z(0,0,0,u)=k_2(u), \; h'_z(0,0,0,u)=k_3(u),\\
f'_{\zeta}(0,0,0,u)=m_1(u),  \; g'_{\zeta}(0,0,0,u)=m_2(u), \; h'_{\zeta}(0,0,0,u)=m_3(u),\\
f'_{\eta}(0,0,0,u)=n_1(u),  \; g'_{\eta}(0,0,0,u)=n_2(u), \; h'_{\eta}(0,0,0,u)=n_3(u).
\end{eqnarray*}
Подставляем полученное значение $e$ в $L$  (обозначение $L$  сохраняем).  Подставим $ \bar{z}= \bar{\zeta}=\bar{\eta}=0$    в $ L_{\bar{z}}, \; L_{\bar{\zeta}}, \;L_{\bar{\eta}}$,  получаем
\begin{eqnarray} \label{dc}
\\
\nonumber 2\,z^2\,\bar{k}_2(u) + 2\,z^3\,\bar{k}_3(u) + 10\,\eta^4\, s_8 \, h(z, \zeta, \eta, u)+24 \,   \bar{c}(u)\,\zeta^2\,z^2+\\
\nonumber 2\,z^5\,\bar{s}_4 \,\bar{k}_3(u)+ 8\,h(z, \zeta, \eta, u)\,\eta^3\, \bar{r}_4 +
4\,z^4\, \bar{r}_2 \,\bar{k}_3(u)+ 2\,\bar{c}(u)\,z^4\,\bar{s}_3+\\
\nonumber 4\,z^4\,\bar{r}_1\,\bar{k}_2(u)-2\,a(u)+2\,z\,\bar{k}_1(u)+
8\,\zeta^2\,\bar{b}(u)+
4\,\zeta\, \bar{a}(u)-\\
\nonumber 2\,d'(u)\,z+2\,f(z, \zeta, \eta, u)+12\,\zeta\,z^2\,\bar{c}(u)-4\,i\,z^2\,a'(u)-\\
\nonumber (4\,i)\,z^4\,\bar{c}'(u)-4\,i\,z^3\,\bar{b}'(u)+
4\,h(z, \zeta, \eta,u)\,\zeta^2\,\eta\,r_3+16\, \bar{c}(u)\,\zeta\,z^3\, \bar{r}_2-\\
\nonumber 8\,i \,z^5\,\bar{r}_2\, \bar{c}'(u)-
8\,i\,\bar{r}_1 \,z^5\, \bar{b}'(u)-
4\,i \, \bar{s}_4 \,z^6\, \bar{c}'(u)+4\,\zeta^2\,\eta\, r_ 3 \, g(z, \zeta, \eta, u)=0\\
\nonumber \\
\nonumber 16\,h(z, \zeta, \eta, u)\,\eta^3\,z\, r_4- 4\,i \,z^5\,\bar{c}'(u)- 2\,i \,z^7\,\bar{s}_4\,c'(u)-4\,i\,z^6\, \bar{r}_1 \,\bar{b}'(u)+\\
\nonumber 2\,z^2\,\bar{m}_2(u)+2\,z^3\,\bar{m}_3(u)+
2\,z\,\bar{m}_1(u)-2\,z^2\,e'(u)+\\
\nonumber 4\,z\,f(z, \zeta, \eta, u)+8\,h(z, \zeta, \eta, u)\,\zeta^2\,\eta\,z\,r_3+2\, \bar{c}(u)\,z^4\,\bar{s}_6+24\, \bar{c}(u)\,\zeta\,z^3 -\\
\nonumber  4\,i \,z^4\,\bar{b}'(u)-
 4\,i \,z^3 \, a'(u)-8 \,i \, z^6 \,\bar{r}_2 \, \bar{c}'(u)+4\,z^4\,\bar{r}_2 \,\bar{m}_3(u)+2\,z^5\,\bar{s}_4\,\bar{m}_3(u)+\\
\nonumber 8\,\zeta\,z\,\bar{a}(u)+16\,\zeta\,z\,\bar{b}(u)+4\,z^4\,\bar{r}_1\,\bar{m}_2(u)=0\\
\nonumber \\
\nonumber 4\,z^4\, \bar{r}_1\, \bar{n}_2(u)+4\,z^4\,\bar{s}_7 \,\bar{c}(u)+2\,z^5\,\bar{s}_4\,\bar{n}_3(u)+6\,z^2\,f(z, \zeta, \eta, u)+\\
\nonumber 4\,z^4\,\bar{r}_2\,\bar{n}_3(u)-4\,i \,z^5\,\bar{b}'(u)-
4\,i\,z^4\,a'(u)-4\,i \,z^7\,\bar{r}_2 \, \bar{c}(u)-\\
\nonumber 2\,z^3\,e'(u)-4\,i\,z^7\,\bar{r}_1 \,\bar{b}'(u)-2\,i\,z^8\,\bar{s}_4\,\bar{c}'(u) - 4\,i \,z^6\, \bar{c}'(u)+\\
\nonumber 2\,z\,\bar{n}_1(u)+2\,z^2\,\bar{n}_2(u)+2\,z^3\,\bar{n}_3(u)=0.
\end{eqnarray}
Из третьего соотношения (\ref{dc}) следует, что $n_1(u)=0$ и что
\begin{eqnarray} \label{f}
3\,f(z, \zeta, \eta, u)=
- \bar{n}_2(u)+z \,(e'(u)-\bar{n}_3(u)) + \\
\nonumber 2\,z^2\, (-\bar{r}_1\, \bar{n}_2(u)-\bar{s}_7 \,\bar{c}(u)-\bar{r}_2\,\bar{n}_3(u)+i\,\bar{a}'(u))+\\
\nonumber z^3 \, (2\,i \, \bar{b}'(u)-\bar{s}_4\,\bar{n}_3(u))+ 2\,i \,z^4\, \bar{c}'(u)\\
\nonumber +2\,i \,z^5\,(\bar{r}_2 \, \bar{c}'(u)+\bar{r}_1 \,\bar{b}'(u))+i\,z^6\,\bar{s}_4\,\bar{c}'(u).
\end{eqnarray}
Подставляя полученное значение $f$ в  первое и второе из соотношений (\ref{dc}) и выделяя старшую компоненту по $\eta$,  получаем, что
\begin{eqnarray}\label{h}
\nonumber    2 \, r_3 \, \zeta^2 \, g(z,\zeta,\eta,u) + (  2 \, r_3 \, \zeta^2 +  5 \,s_8 \, \eta^2)\, h(z,\zeta,\eta,u) =0 , \\
 (r_3 \, \zeta^2 + 2 \, r_4 \, \eta^2) \, h(z,\zeta,\eta,u)=0.
\end{eqnarray}
Оставшиеся после выполнения (\ref{h})  условия, обеспечивающие выполнение  (\ref{dc}),  сводятся к  следующей системе соотношений
\begin{eqnarray} \label{sdc2}
a=b=c=n_1=n_2=k_3=0,   \quad \quad\quad\quad\quad\quad\quad\quad  \\
\nonumber
d'=-2\,{\rm Re} \, n_3,    \quad  k_2=\frac{2}{3} \, r_2 \, n_3,   \quad  m_2=\frac{1}{3} \, (n_3-\bar{n}_3),  \quad m_3 =\frac{4}{3} \, r_2 \, n_3,  \\
\nonumber
r_1 \, r_2 \, n_3 = 0, \quad   s_4 \, n_3=0,   \quad  (r_1-s_4+4 \, r_2^2) \, n_3 = r_1 \, \bar{n}_3.\quad\quad\quad
\end{eqnarray}

Отметим, что при этом
$$e(z,\zeta,\eta,u)=d(u), \quad f(z,\zeta,\eta,u) = \frac{d'(u)-\bar{n}_3(u)}{3} \, z -\frac{2}{3} \, \bar{r}_2 \, \bar{n}_3(u) \, z^2 -\frac{\bar{s}_4 \, \bar{n}_3(u)}{3} \, z^4.$$

Рассмотрим несколько случаев.\\
{\it Случай 1:}  Пусть $r_3 \neq 0$,  тогда из (\ref{h}) сразу следует, что $g=h=0$. Оставшиеся соотношения позволяют заключить, что $f=0$ и что $d$ - это вещественная константа.  Итак, в этом случае ${\rm Ker} \,L=V^0=\{ (0,0,0,d_0)\}$, причём $d_0 \in \mathbf{R}$.   \\
{\it Случай 2:}  Пусть $r_3 = 0$.
 Пусть $(r_4,s_8) \neq 0$,   ({\it случай 2.1}), тогда из (\ref{h}) следует, что $h=0$. Возвращаясь к соотношениям (\ref{dc}), получаем
 $n_3=n_2=d'=0$.  Откуда $f=0,$ а $e$ - вещественная константа. Обозначим $g''_{zz}(0,0,0,u)$ через $k_{22}(u)$.
 Вычисляя теперь $L''_{\bar{z} \bar{z}}$ при $\bar{z}=\bar{\zeta}= \bar{\eta}=0$, получаем
 $$
 g(z,\zeta,\eta,u) = 2 \,i \,z^5\,\bar{r}_1 \, \bar{k}'_2(u)-z^4 \, \bar{r}_1 \, \bar{k}_{22}(u)+i\,z^3\, \bar{k}'_2(u)  -(1/2) z^2 \, \bar{k}_{22}(u)-4 \, \zeta^2 \, \bar{k}_2(u).$$
 Подставляя это значение $g$ в $L$ и приравнивая к нулю коэффциенты при  $z^2 \, \bar{\zeta}^2$ и $z^2 \, \bar{z} \, |\zeta|^2$, получаем,  что $k_2=k_{22}=0$, т.е. $g=0$.  Т.е. ${\rm Ker} \,L=V^0$.

 Пусть $r_4=s_8= 0$,   ({\it случай 2.2}).   Обозначим
 $$h''_{zz}(0,0,0,u)=k_{32}(u), \; h'''_{zzz}(0,0,0,u)=k_{33}(u).$$
  Вычисляя
 $L''_{\bar{z} \bar{z}}$ при $\bar{z}=\bar{\zeta}= \bar{\eta}=0$,  как и ранее, получаем выражение для $g$,  вычисляя
  $L'''_{\bar{z} \bar{z} \bar{z}}$ при $\bar{z}=\bar{\zeta}= \bar{\eta}=0$,  получаем выражение для $h$. После чего анализ
  младших коэффициентов $L$ даёт  $n_3=d'=k_2=k_{22}=k_{32}=k_{33}=0$,  откуда получаем $f=g=h=0$,  $e$ - вещественная константа.  Т.е. ${\rm Ker} \,L=V^0$.

  Таким образом нами доказана следующая лемма.

 \vspace{2ex}

 {\bf Лемма 4:}  Если $L(\phi)=0$, то $\phi = (0,0,0,d_0)$, где $d_0$ - вещественная постоянная. В частности, на пространстве $V_5$ ядро  тривиально.

\vspace{2ex}

Выясним, как устроены младшие струи отображения  $\Gamma$ на  $\tilde{\Gamma}$, сохраняющего ноль на месте.
Выделяя компоненты (\ref{mid}) веса один и два, сразу получаем, что
$$ e_1=0, \quad e_2=|\lambda|^2 \, w, \quad f_1 = \lambda \, z, \quad   \lambda \in \mathbf{C}^{*}.$$
Пусть, далее,
$$e_3=|\lambda|^2 \, (d_3+d_1 \, w),  \quad f_2=\lambda \, (a \, w + a_2 ),   \quad g_1=b_1,$$
где $d_3,d_1,a_2,b_1$ - однородные голоморфные формы от $(z,\zeta,\eta)$ соответствующих степеней, $a$ - постоянная.
Третья весовая компонента (\ref{mid}) имеет вид
\begin{eqnarray*}
|\lambda|^2 \, {\rm Im} [ i \, 2 \, {\rm Re}(z^2 \, \bar{\zeta} + d_1(z,\zeta,\eta) \, (u+i \, |z|^2))] = \\
|\lambda|^2 \, 2 \, {\rm Re} [(a \, (u+i \, |z|^2)+a_2(z,\zeta,\eta)) \, \bar{z}] + 2 \, {\rm Re} [\lambda^2 \, z^2 \, \bar{b}_1(\bar{z},\bar{\zeta},\bar{\eta})].
\end{eqnarray*}
Выделяя члены, линейные по $u$, получаем $d_1=2 \, i \, \bar{a} \, z$. Теперь среди членов бистепени$(2,1)$ выпишем по отдельности компоненты, линейные по $\bar{z}$, по  $\bar{\zeta}$ и по $\bar{\eta}$. Получим
\begin{eqnarray*}
a_2(z,\zeta,\eta)=(2 \, i \, \bar{a} - \frac{\lambda}{\bar{\lambda}}\, \bar{b}_1^1) \, z^2, \quad
b_1^2=\frac{\lambda}{\bar{\lambda}},  \quad  b_1^3=0,
\end{eqnarray*}
где $b_1= b_1^1 \, z + b_1^2 \, \zeta + b_1^3 \,\eta$. Здесь и далее верхними индексами указываем на связь коэффициентов и переменных.  Пусть  $b_1^1=\alpha \, \lambda / \bar{ \lambda}$,  получаем
\begin{eqnarray*}
e_3=2 \, i \, |\lambda|^2 \, \bar{a} \, z \, w, \quad f_2= \lambda \, (a \, w + (2 \, i \, \bar{a}- \bar{\alpha}) \, z^2), \quad
g_1=\frac{\lambda}{\bar{\lambda}} (\alpha \, z + \zeta).
\end{eqnarray*}

Пусть, далее,
\begin{eqnarray*}
e_4=|\lambda|^2 \,(d_4+d_2 \, w + d_0 \, w^2),   \quad f_3=\lambda \, (a_3+a_1 \, w) , \\
 g_2 = \frac{\lambda}{\bar{\lambda}}(b_2+b_0 \, w), \quad h_1=   \frac{\lambda}{\bar{\lambda}^2} \, c_1,
\end{eqnarray*}
где коэффициенты -- это голоморфные однородные формы от  $(z,\zeta,\eta)$ соответствующих степеней. Выпишем компоненту (\ref{mid}) веса четыре (общий множитель $|\lambda|^2$  убираем).
\begin{eqnarray} \label{c4}
\nonumber {\rm Im} [d_4(z,\zeta,\eta)+d_2(z,\zeta,\eta) \, (u+i \,|z|^2 ) + d_0 \, (u^2+2 \, i \,|z|^2 \, u -|z|^4 )+\\
\nonumber  -   2 \,  \bar{a} \, z \, (z^2 \, \bar{\zeta}+ \bar{z}^2 \, \zeta) +  \, i \, ( z^3 \, \bar{\eta}+\bar{z}^3 \, \eta)+4 \,|z|^2 \, |\zeta|^2) ]=\\
  2 \, {\rm Re}[(a_3(z,\zeta,\eta)  + a_1(z,\zeta,\eta) \, (u+i \,|z|^2 )) \, \bar{z}]  +\\
 \nonumber |a|^2 \, (u^2+|z|^4 )- 2 \, {\rm Re}[(a \, (u+i \,|z|^2 ) \, (2 \, i \, a+ \alpha) \,\bar{ z}^2)]+ |2 \, i \ a + \alpha|^2 \, |z|^4+ \\
  2 \, {\rm Re}[ a \,(2 \,i \,  {\rm Re}(z^2 \, \bar{\zeta}) )] +2 \, {\rm Re}[2 \,( a \, (u+i \, |z|^2) + (2 \, i \, \bar{a}- \bar{\alpha}) \, \nonumber  z^2) \, z \, ( \bar{\alpha} \, \bar{z} + \bar{\zeta})] +\\
 \nonumber  2 \, {\rm Re}[(b_2(z,\zeta,\eta)+b_0 \, (u+i \, |z|^2)) \, \bar{z}^2]+
2 \, {\rm Re}[c_1(z,\zeta,\eta) \, \bar{z}^3] + 4 \, |z|^2 \, |\alpha \, z + \zeta|^2.
\end{eqnarray}
Отделяя в (\ref{c4}) коэффициент при $u^2$, получаем  ${\rm Im} \, d_0=|a|^2$. Положим $d_0 = \gamma + i \, |a|^2$.
Отделяя в (\ref{c4}) коэффициент при $u$, получаем
\begin{eqnarray*}
{\rm Im}[ d_2(z,\zeta,\eta)] +|a|^2 \, |z|^2 =  2 \, {\rm Re} [a_1(z,\zeta,\eta) \, \bar{z}] - \\
2 \, {\rm Re} [a \, (2 \, i \, a + \alpha) \, \bar{z}^2 ] +  2 \, {\rm Re} [a \, z \, (\bar{\alpha} \, \bar{z} + \bar{\zeta})] + 2 \, {\rm Re} [b_0 \, \bar{z}^2]
\end{eqnarray*}
Откуда получаем
\begin{eqnarray*}
d_2=(2 \, i \, \bar{a}^2 - \bar{a} \, \bar{\alpha}+\bar{\beta}) \, z^2, \\
a_1 =  ((|a|^2- 2 \, {\rm Re}(a \,  \bar{\alpha}))+i \,\delta) \, z +  \bar{a} \, \zeta,
\end{eqnarray*}
где $\beta = b_0, \; \delta ={\rm Im} \, a_1^1$.  Компонента бистепени $(4,0)$ сразу дает $d_4=0$. В бистепени $(3,1)$
получаем
\begin{eqnarray*}
i \, d_2(z,\zeta,\eta) \, |z|^2 - \bar{a} \, z^3 \, \bar{\zeta} + z^3 \, \bar{\zeta} =
 a_3(z,\zeta,\eta) \, \bar{z} +\\
 i \,  \bar{a} \,( - 2 \, i \, \bar{a} + \bar{\alpha}) z^3 \, \bar{z} - i \,  \bar{\beta} \, z^3 \, \bar{z}+ z^3 \, \overline{c_1(z,\zeta,\eta)}.
\end{eqnarray*}
Откуда получаем
\begin{eqnarray*}
a_3= (- 4 \, \bar{a}^2 - 2 \, i \, \bar{a} \bar{\alpha} + 2 \, \bar{\beta}- \bar{\nu}) \, z^3,\\
c_1 = \nu \, z - a \, \zeta +\eta,
\end{eqnarray*}
где $\nu = c_1^1$. В бистепени $(2,2)$ имеем
\begin{eqnarray*}
2 \, {\rm Re} [ - i \, \bar{a} z \, \zeta \, \bar{z}^2 + i \, a_1 \, z \,  \bar{z}^2 + 2 \, i \, a \, z^2 \, \bar{z} \, (\bar{\alpha}\, \bar{z}+\bar{\zeta}) + b_2 \, \bar{z}^2 + 4 \, \alpha \, z^2 \,\bar{z} \, \bar{\zeta}]+\\
(|a|^2 + 4 \, |\alpha|^2 +|2 \, i \, a + \alpha|^2) \, |z|^4=0.
\end{eqnarray*}
Откуда получаем
\begin{eqnarray*}
b_2= (\delta- 2 \, |\alpha|^2 -\frac{|a|^2}{2}-\frac{1}{2}\, |2 \, i \, a + \alpha|^2- 2 \, {\rm Re}(i \, a \, \bar{\alpha} + i  \, \kappa)) \, z^2 + (i \, a + 2 \alpha) \, z \, \zeta, \\
c_1= \nu \, z - a \,\zeta +\eta,
\end{eqnarray*}
где $\kappa = {\rm Im} \, b_2^{11}$.

Итак, нами доказана следующая лемма.

\vspace{3ex}

 {\bf Лемма 5:}(a) Всякое локально обратимое отображение $\Gamma$ на  $\tilde{\Gamma}$, сохраняющее начало координат, представимо в виде композиции отображения вида
\begin{eqnarray*}
z \rightarrow \lambda \,(  z + a \, w +
(2 \, i \, \bar{a} - \bar{\alpha})\, z^2 +(-4\,\bar{a}^2+2\,i\,\bar{a}\,\bar{\alpha}+2\,\bar{\beta}-\bar{\nu})\, z^3+\\
((|a|^2-2 \, {\rm Re}(\,a\, \bar{\alpha})+i \, \delta) \, z + \bar{a} \, \zeta)  \, w),  \qquad\qquad\qquad\qquad\\
\zeta  \rightarrow \frac{\lambda}{\bar{\lambda}} \, (\zeta  + \alpha \, z + \tau \, z^2 + (i \, a + 2 \alpha) \, z \, \zeta + \beta \, w),\qquad\qquad\\
\eta \rightarrow \frac{\lambda}{\bar{\lambda}^2} \, (\eta +\nu \, z - a \,\zeta), \qquad\qquad\qquad\qquad\qquad\qquad\qquad\\
w   \rightarrow |\lambda|^2 \, (w + 2 \, i  \, \bar{a} \, z \, w + (2 \, i \, \bar{a}^2 - \bar{a} \, \bar{\alpha}+\bar{\beta}) \, z^2 \, w + ( \gamma + i \, |a|^2) \, w^2), \\
\mbox{ где } \quad \tau= (\delta- 2 \, |\alpha|^2 -\frac{|a|^2}{2}-\frac{1}{2}\, |2 \, i \, a + \alpha|^2- 2 \, {\rm Re}(i \, a \, \bar{\alpha} + i  \, \kappa))
\end{eqnarray*}
и отображения
$$z   \rightarrow  z +O(4), \;  \zeta  \rightarrow \zeta + O(3) , \;  \eta \rightarrow \eta +O(2),  \; w   \rightarrow w +O(5).$$
(b) Причём
$$\lambda \in \mathbf{C}^{*}, \quad  a, \; \alpha, \; \beta, \; \nu  \in \mathbf{C},   \quad   \gamma, \; \delta, \; \kappa \in \mathbf{R}.$$
Что дает  13 вещественных параметров.\\
(c)  Такое отображение однозначно определяется заданием 3-струи в начале координат.
\vspace{3ex}

Теперь мы готовы доказать следующее утверждение.

\vspace{3ex}

{\bf Утверждение 6:}  Если $\Gamma$ --  вещественно аналитическая гиперповерхность пространства $ \mathbf{C}^4$, которая в точке общего положения является 3-невырожденной, то размерность группы локальных голоморфных автоморфизмов {\it в любой точке} не превосходит 20.\\
{\it Доказательство:}  Размерность группы в призвольной точке не превосходит суммы  размерности гиперповерхности и размерности стабилизатора в точке общего положения.   Размерность гиперповерхности равна 7. Размерность стабилизатора, в силу теоремы  1 и лемм 2,3,4 и 5 не превосходит 13.  Поскольку $ 7 + 13 =20$,   утверждение доказано.\\

\vspace{7ex}

{\bf  Общие 2-невырожденные гиперповерхности}

\vspace{3ex}

Обозначим координаты в $\mathbf{C}^4$ через $(z=(z_1,z_2), \zeta, w=u+i \, v)$ и перейдем к рассмотрению {\it 2-невырожденного} случая.Так же как и выше,  мы можем ограничиться "жесткими"\  гиперповерхностями.
Пусть  $\Gamma$ --  равномерно 2-невырождена в окрестности ${\xi} \in  \Gamma$.
Форма Леви равномерно 2-невырожденной гиперповерхности повсюду имеет имеет минимальное вырождение, а именно, ее ранг ранг равен 2.   Таким образом  мы можем записать локальное уравнение $ \Gamma$ в виде
\begin{equation}\label{e}
  v = < z, \bar{z} >  +  F_3(z,\bar{z},\zeta, \bar{\zeta})+F_4 (z,\bar{z},\zeta, \bar{\zeta})+ \dots
\end{equation}
где $F_j$ -- однородный вещественный полином степени $j$,  а $ < z, \bar{z} >$ -- невырожденная эрмитова форма от переменного $z \in \mathbf{C}^2$.  Простыми треугольно-полиномиальными заменами переменных $z$ и $w$ можно добиться того, что правая часть уравнения  $\Gamma$
$$F=< z, \bar{z} >  + F_3+F_4+\dots$$
 не будет содержать плюригармонических слагаемых (т.е. слагаемых бистепеней $(m,0)$ и $(0,m)$) и слагаемых, линейно зависящих от $z$ и $\bar{z}$, за исключением формы  $< z, \bar{z} >$.  Выпишем слагаемые, которые после этого останутся в $F_3$ и $F_4$. Имеем
\begin{eqnarray}  \label{34}
 F_3=2 \, {\rm Re} ( K(z,z) \, \bar{\zeta} + A_1(z) \, |\zeta|^2 + A_2 \, \zeta^2 \, \bar{\zeta}),  \quad\quad\quad\quad\quad \quad\quad \\
\nonumber  F_4= 2 \, {\rm Re} ( (P(z,z,\bar{z})  + Q(z,z,z)) \, \bar{\zeta}+R(z,z) \bar{\zeta}^2) + S(z,\bar{z}) |\zeta|^2 +T(z,z,\bar{z},\bar{z})+\\
 \nonumber 2 \, {\rm Re} ( B_1(z,z) \, |\zeta|^2   +  B_2(z) \, \zeta^2 \, \bar{\zeta} +  B_3(z) \, \zeta  \, \bar{\zeta}^2 ). \quad\quad\quad\quad\quad\quad\quad
 \end{eqnarray}

 \vspace{3ex}

 Отметим,   что из условия 2-невырожденности $\Gamma$ следует, что форма $K(z,z)$ не равна нулю тождественно.
Для  дальнейших вычислений нам потребуется приведение пары форм $(< z, \bar{z} >, K(z,z))$  на $\mathbf{C}^2$ комплексно линейными заменами к  виду, содержащему минимум параметров.   Имеет место следующая классификация.

 \vspace{3ex}

 {\bf Лемма 7:}  Пусть $< z, \bar{z} >$   невырождена, а $K(z,z)$ отлична от тождественного нуля, тогда невырожденным комплексно линейным преобразованием можно привести эту пару к одной из следующего списка:\\
 \begin{eqnarray*}
 \nonumber
 (1.)  \quad   (|z_1|^2+|z_2|^2, \; k \, z_1^2+ m \, z_2^2) ,  \quad  k, \, m   >0, \; k \neq m,\\
 (2.)  \quad   (|z_1|^2+|z_2|^2, \; k \, (z_1^2  + z_2^2)) ,  \quad  k >0, \\
 (3.)  \quad   (|z_1|^2+|z_2|^2, \; k \, z_1^2) ,  \quad  k >0,  \quad \quad \quad  \quad \\
 (4.)  \quad   (|z_1|^2-|z_2|^2, \; k \, z_1^2+ m \, z_2^2) ,  \quad  k, \, m >0, \; k \neq m,\\
 (5.)  \quad   (|z_1|^2-|z_2|^2, \; k \,( z_1^2 +  z_2^2)) ,  \quad  k >0, \quad \quad \quad  \quad \\
 (6.)  \quad   (|z_1|^2-|z_2|^2, \; k \, z_1^2) ,  \quad  k  >0,  \quad \quad \quad \quad\\
 (7.)  \quad   ( 2 \, {\rm Re} \,(z_1 \, \bar{z}_2) ,\;  z_1^2+ m \, z_2^2) ,  \quad   m \notin \mathbf{R},  \quad \quad \\
(8.)  \quad   ( 2 \, {\rm Re} \,(z_1 \, \bar{z}_2) ,\;  z_1^2+ m \, z_2^2) ,  \quad   m \in \mathbf{R}^{*},  \quad \quad \\
(9.)  \quad   ( 2 \, {\rm Re} \,(z_1 \, \bar{z}_2), \;  z_1^2 ) . \quad \quad \quad \quad \quad \quad  \quad \quad \quad \\
 \end{eqnarray*}
 {\it Доказательство:} Пусть $< z, \bar{z} >$  положительно определена и $\nu$ -- собственный вектор оператора, заданного матрицей
 $$ \left[ \begin {array}{cc} k&l\\ \noalign{\medskip}l&m\end {array}\right].$$
 Выбирая в качестве первого вектора нового базиса вектор
 $$\frac{\nu}{\sqrt{< \nu, \bar{\nu} >}}$$
 и, подбирая второй из условия ортонормированности, получаем в зависимости от ранга $K$ пары (1.), (2.) и (3.).  Положительности параметров  $k$ и $m$ можно добиться поворотами в плоскостях $z_1$ и $z_2$.\\
 Пусть $< z, \bar{z} >$ -- имеет сигнатуру $(1,1)$.   Если оператор имеет собственный вектор $\nu$, т.ч.  $< \nu, \bar{\nu} > \neq 0$, то годится то же самое рассуждение и это дает пары (4.), (5.) и (6.).\\
Пусть $(e_1,e_2)$  --  базис $\mathbf{C}^2$, в котором $K(z,z)$ -- диагональна, т.е. $K(z,z)= k \, z_1^2+ m \, z_2^2$ причем $< e_1, \bar{e}_1 > =< e_2, \bar{e}_2 > =0$. Тогда если $z= z_1 \, e_1 + z_2 \, e_2$, то
$$ <z, \bar{z}>= 2 \, {\rm Re} (<e_1, \bar{e}_2> \, z_1 \, \bar{z}_2).$$
После растяжения по $z_1$ эрмитова форма принимает вид $<z,\bar{z}>= 2 \, {\rm Re} ( z_1 \, \bar{z}_2)$.
Используя преобразование
$$ z_1 \rightarrow \lambda \, z_1 , \quad z_2 \rightarrow \frac{z_2}{\bar{\lambda }},  \quad \lambda \in \mathbf{C}^{*},$$
которое не меняет эрмитовой формы получаем пары (7.), (8.) и (9.).  Лемма доказана.

\vspace{3ex}

{\bf Лемма 8:}    Если форма Леви гиперповерхности $\Gamma$, заданной уравнением (\ref{e}), тождественно вырождена, то $F_3$ и $F_4$ можно записать в виде (\ref{34}), причем    $A_1=A_2=B_1=B_2=B_3=0,$ а
форма $S$ в зависимости от номера пары из леммы 2 имеет следующий вид:
\begin{eqnarray*}
(1.) \; S=4 \, (k^2 \, |z_1|^2+m^2 \, |z_2|^2), \;  (2.) \; S=4 \, k^2 \, (|z_1|^2+ |z_2|^2), \;  (3.) \;  S=4 \, k^2 \, |z_1|^2, \qquad \qquad\qquad\qquad\qquad\\
(4.) \;  S=4 \, (k^2 \, |z_1|^2 - m^2 \, |z_2|^2), \; (5.) \;  S=4 \, k^2 \, (|z_1|^2 - |z_2|^2), \; (6.) \;  S=4 \, k^2 \, |z_1|^2, \quad\qquad\qquad\qquad\qquad\\
(7.) \;  S=4 \, (\bar{m} \, z_1 \, \bar{z}_2+m \, z_2 \, \bar{z}_1), \; (8.) \;  S=4 \,m \, ( z_1 \, \bar{z}_2+ z_2 \, \bar{z}_1), \; (9.) \;  S=0. \qquad\qquad\qquad \qquad\qquad\qquad
\end{eqnarray*}
{\it Доказательство:}  Вычисляя определитель матрицы комплексного гессиана по переменным $(z_1,z_2,\zeta)$  и отделяя в нем компоненты степени один, получаем, что  $A_1=A_2=0$ .  Отделяя, далее, компоненты степени два, получаем, что  $B_1=B_2=B_3=0$ и указанный вид формы $S$.  Лемма доказана.

  \vspace{3ex}

 Теперь можем написать
\begin{eqnarray*}
\nonumber F_3=   2 \, {\rm Re}[ K(z,z)  \, \bar{\zeta} ] ,\qquad \qquad\qquad\qquad\qquad\qquad \qquad\qquad\qquad\qquad\\
\nonumber  F_4= 2 \, {\rm Re} [ (P(z,z,\bar{z})  + Q(z,z,z)) \, \bar{\zeta}+R(z,z) \bar{\zeta}^2] + S(z,\bar{z}) |\zeta|^2 +T(z,z,\bar{z},\bar{z}).
 \end{eqnarray*}
т.е. уравнение гиперповерхности имеет вид
\begin{eqnarray}\label{G2}
\nonumber
 v=< z, \bar{z} >+ 2 \, {\rm Re}[ K(z,z)  \, \bar{\zeta} ]+\\
 2 \, {\rm Re} [ (P(z,z,\bar{z})  + Q(z,z,z)) \, \bar{\zeta}+R(z,z) \bar{\zeta}^2] + S(z,\bar{z}) |\zeta|^2 +T(z,z,\bar{z},\bar{z}) +O(5).
\end{eqnarray}
\vspace{3ex}

Введем в пространстве степенных рядов от $(z,\bar{z},\zeta, \bar{\zeta},u)$, а также от $(z,\bar{z},\zeta,\bar{\zeta}, w,\bar{w})$, градуировку, назначая веса переменным
$$ [z]=[\bar{z}]=[\zeta]=[ \bar{\zeta}]=1, \; [w]=[\bar{w}]=[u]=2.$$
Пусть $\Gamma$ и $\tilde{\Gamma}$ -- гиперповерхности, заданные уравнениями
\begin{eqnarray} \label{hs}
 v = < z, \bar{z} >  + 2 \, {\rm Re} \,( K(z,z)  \, \bar{\zeta} )  +O(4),\\
 \nonumber
v = < z, \bar{z} >  + 2 \, {\rm Re} \, (\tilde{K}(z,z)  \, \bar{\zeta} )  +O(4)
\end{eqnarray}
и
$$\phi=( \; z \rightarrow f=f_1 + f_2 + O(3), \;  \zeta \rightarrow g=g_1 + O(2), \;  w \rightarrow h=h_1 + h_2 + h_3+O(4) \;)$$
-- локально обратимое голоморфное отображение первой на вторую, оставляющее начало координат на месте. Причем компоненты
координат отображения -- это компоненты фиксированного веса и $O(j)$ --  это сумма слагаемых веса не ниже $j$.

То, что это отображение переводит $\Gamma$ в $\tilde{\Gamma}$, аналитически можно записать в виде следующего соотношения
\begin{eqnarray} \label{me}
 \nonumber
\Theta = -2 \,  {\rm Im} \, h  +2 \, <f,\bar{f}> + 4 \, {\rm Re} \,( \tilde{K}(f,f)  \, \bar{g} )+ \tilde{F}_4+ \dots  \\
 \nonumber    \mbox{   при     }  w=u+i \,(< z, \bar{z} >  + 2 \, {\rm Re} \, ( K(z,z)  \, \bar{\zeta} )  +F_4+O(5)).\\
  \end{eqnarray}
Отделяя в этом соотношении компоненту веса 1, получаем $h_1=A(z)+B \, \zeta=0$. \\
Пусть $h_2= \Phi_2(z,\zeta)+ \rho \, w,   \;  f_1= C \, z + d \, \zeta$, где $\Phi_2$ -- форма степени два от $(z,\zeta)$.
Отделяя в (\ref{me}) компоненту веса 2, получаем $\Phi_2(z,\zeta)=0,   \quad  <C \,z, \overline{C \, z}>= \rho \, < z, \bar{z} > $.
Т.е. $f_1=C \,z,     \quad h_2 = \rho \, w$.   Отметим, что в силу обратимости отображения матрица $C$ -- невырождена и $\rho \neq 0$.\\
Пусть, далее,
\begin{eqnarray*}
h_3= \rho (\Phi_3 + (A(z)+B \, \zeta) \, w) ,  \quad f_2= C(a \, w + b(z,z) + c(z) \, \zeta +d \, \zeta^2), \\
 g_1 = <z, \bar{\alpha}> + \beta \, \zeta,
\end{eqnarray*}
где $\Phi_3$ -- форма степени три от $(z,\zeta)$. Тогда, отделяя  в (\ref{me}) компоненту веса 3, получаем
\begin{eqnarray*}
h_3=2 \, i \, \rho \,   <z,\bar{a}> \,w, \quad  f_2=C(aw+2 \, i\,  <z,\bar{a}> \, z -  K(z,z) \, \mu), \\
g_1=  <z, \bar{\alpha}> + \beta \, \zeta,
\end{eqnarray*}
причем в силу обратимости отображения $\beta \neq 0$.

Из наших вычислений видно, что при определении весовой $j$-струи удобно принять следующую точку зрения.
$$\phi = \sum \phi_j, \quad \phi_j=(f_{j-1},  \, g_{j-2}, \, h_j)$$
Т.е. весовая $j$-струя отображения понимается как набор струй координат, где у $h$ берем $j$-ю весовую cтрую,  у $f$ -- $(j-1)$-ю, у $g$ берем $(j-2)$-ю.  Проведенное выше вычисление дает описание действия голоморфных отображений на 3-струю уравнения  гиперповерхности  вида (\ref{hs}).

\vspace{3ex}

{\bf Лемма 9:}\\
(a) На совокупности весовых 3-струй 2-невырожденных гиперповерхностей вида (\ref{hs}) псевдогруппа локально обратимых голоморфных отображений, сохраняющих начало координат, действует следующим образом:
\begin{eqnarray*}
\nonumber
z \rightarrow  C(z+ aw+2 \, i  \, <z,\bar{a}> \, z -  K(z,z) \, \alpha) + O(3), \\
\nonumber
 \zeta \rightarrow  <z, \bar{\alpha}> + \beta \zeta   +O(2),\\
\nonumber
w \rightarrow   \rho \, ( w + 2 \, i \,  <z,\bar{a}> \,w)  +O(4),
\end{eqnarray*}
причем $C \in {\rm GL}(2,\mathbf{C}), \; \rho \in \mathbf{R}^{*}, \;  a, \; \alpha \in \mathbf{C}^2, \; \beta \in \mathbf{C}^{*} $,  а также
\begin{eqnarray} \label{r}
  <z, \bar{z}>=\rho \, <C^{-1} z, \overline{C^{-1} z}>,  \quad  \tilde{K}(z,z) = \frac{\rho}{\bar{\beta}}K(C^{-1} \,z,C^{-1} \,z)
\end{eqnarray}
(b)  Любое обратимое голоморфное отображение $\Gamma$ на $\tilde{\Gamma}$, сохраняющее начало координат,   с точностью до этого действия имеет вид
\begin{equation}\label{nm}
z \rightarrow z + O(3), \quad  \zeta \rightarrow \zeta + O(2), \quad  w \rightarrow w +O(4).
\end{equation}

\vspace{1ex}

Пусть речь идет об отображении гиперповерхности $\Gamma$ на себя. Тогда рассмотрим подгруппу группы автоморфизмов  $\Gamma$,  состоящую из линейных автоморфизмов вида
\begin{eqnarray}\label{la}
\nonumber    G_0 = \{ (z  \rightarrow C \, z,  \;\zeta  \rightarrow \beta \, \zeta,  \;w  \rightarrow \rho \, w)\}  \mbox{  с условием}\\
 <C \, z, \overline{C\, z}>=\rho \, <z, \bar{z}> ,  \quad  K(C \, z,C \, z) = \frac{\rho}{\bar{\beta}}K(z,z).
\end{eqnarray}
Вычислим размерность этой группы для каждой из шести пар форм, перечисленных в лемме 8.

\vspace{3ex}

 {\bf Лемма 10:}  Пусть $G_0^j$ -- это группа $G_0$ для $j$-й пары из списка леммы 8. Тогда
 \begin{eqnarray*}
  {\rm dim} \, G_0^1 =2, \;  {\rm dim} \, G_0^2 =3, \; {\rm dim} \, G_0^3 =3, \;  {\rm dim} \, G_0^4 =2,
 \;{\rm dim} \, G_0^5 =3, \\
 {\rm dim} \, G_0^6 =3, {\rm dim} \, G_0^7 =3, \;  {\rm dim} \, G_0^8 =3, \;{\rm dim} \, G_0^9 =3,
 \end{eqnarray*}
 т.е. в любом случае ${\rm dim} \, G_0 \leq 3$.\\
 {\it Доказательство:} {\it Пары (1.),(2.),(3.)}.  $<z,\bar{z}>=|z_1|^2+|z_1|^2$,   тогда $C= \lambda \, U, \; \rho =|\lambda|^2,$ где $U \in SU(2)$, а $\lambda \in \mathbf{C}^{*}$.  Записывая $U$ в виде
 $$ \left[ \begin {array}{cc} p&q\\ \noalign{\medskip}-\bar{q}&\bar{p}\end {array}\right], \quad \mbox{ где         }  |p|^2+|q|^2=1.$$
 и подставляя это во второе соотношение, получаем, что
 $$ (k \, p^2+m \, \bar{q}^2, \; - 4 \, i \, {\rm Im}\, p \, q, \;  k \, q^2 + m \, \bar{p}^2)=\frac{|\lambda|^2}{\bar{\beta}}\, (k,0,m).$$
 Откуда  ${\rm Im} \, p \, q =0$,   т.е. $q= \sigma \, \bar{p}$,  где $\sigma$ вещественно, при этом $|p|^2 \, (1+\sigma^2)=1$. Пусть
 $p=\exp(i \, \phi) /\sqrt{1+\sigma^2}$, тогда имеем
 $$  \exp(4 \,i \, \phi)=\frac{k \, \sigma^2 + m}{k + m \, \sigma^2}\, \frac{k}{m}.$$
 Откуда следует, что для $U$ из малой окрестности единицы  $\phi=0$. Далее имеем: либо $k=m$, либо $\sigma=0$.
 Так же получаем ответ для третьей пары.  Для пары (1.) свободный параметр -- это $\lambda$, для (2.) --
  $(\lambda,\sigma)$, для (3.) -- $(\lambda,\phi)$.

 {\it Пары (4.),(5.),(6.)} Рассматриваются вполне аналогично с учетом того, что $U$ - псевдоунитарная матрица вида
 $$ \left[ \begin {array}{cc} p&q\\ \noalign{\medskip}\bar{q}&\bar{p}\end {array}\right], \quad \mbox{ где          }  |p|^2-|q|^2=1.$$
 Свободные параметры -- те же.

 {\it Пары (7.),(8.),(9.)}. Для этой эрмитовой формы псевдоунитарная матрица с единичным определителем, близкая к единичной, имеет вид
 $$ \left[ \begin {array}{cc} p &{i \,\sigma \, p}\\ \noalign{\medskip}{\frac{i \, r}{(1+r \sigma) \, p}}&{\frac{1}{(1+r \sigma) \, p}}\end {array}\right], \quad \mbox{ где   }  r, \sigma  \in \mathbf{R},   p>0 .$$
Откуда получаем значения размерностей.
Для пары (7.) свободные параметры -- это $(\lambda,p)$, для (8.) --   $(\lambda,r)$, для (9.) -- $(\lambda,p)$.
Лемма доказана.

\vspace{5ex}

Зафиксируем  гиперповерхности  $\Gamma$ и $\tilde{\Gamma}$ вида  (\ref{hs}) и дадим оценку числа параметров, от которых зависит отображение  одной на другую вида (\ref{nm}) в соответствии со схемой рекурсии глубины $k=2$ (теорема 1). С этой целью опишем вид  $\mu$-й компоненты соотношения (\ref{me}).   При этом явно выпишем слагаемые, зависящие от $\phi_{\mu}$ и $\phi_{\mu-1}$, игнорируя члены, зависящие от $\phi_{\nu}$ при $\nu \leq \mu-2$. Пусть $f=(f^1,f^2)$. Введем также следующее обозначение: $\Delta \, \psi(u) =  2\, i\, {\rm Re} (K(z,z)\,\bar{\zeta})  \, \psi'(u)$.

\vspace{3ex}

{\bf Лемма 11:} $\mu$-я весовая компонента выражения (\ref{me}) $\Theta_{\mu}$ имеет вид
 \begin{eqnarray} \label{L2}
\nonumber  \Theta_{\mu} = L_1(\phi_{\mu})+L_2(\phi_{\mu-1}) + \theta_{\mu}(\phi_{\nu < \mu-1}), \qquad  \mbox{ причем } \qquad \qquad \qquad \qquad \qquad\qquad\qquad \\
\nonumber   L_1(\phi) =  2 \, {\rm Re} ( i \, h + 2 \, <f,\bar{z}> + 2 \,  \bar{K}( \bar{z}, \bar{z}) \, g), \qquad   L_2(\phi) = \Delta L_1 (\phi) + \qquad\qquad \qquad\qquad\\
  2 \, {\rm Re} \{ 4 \, K(f,z) \, \bar{\zeta} +
  2\,( \bar{P}( \bar{z}, \bar{z},z)+\ \bar{Q}( \bar{z}, \bar{z},\bar{z}) +
  \nonumber 2 \, \bar{R}( \bar{z}, \bar{z}) \, \zeta + S(z,\bar{z}) \, \bar{\zeta}) \, g \},  \qquad\qquad\qquad\qquad\qquad \\
 \mbox{ где }  \quad  w=u+i \, <z,\bar{z}>.  \qquad\qquad\qquad\qquad\qquad\qquad\qquad\qquad\qquad
  \end{eqnarray}

\vspace{3ex}
Отметим, что выражение $L (\phi)= L_1 (\phi)+ L_2(\phi)$  линейно по $\phi$  и не зависит от $\mu$.
Пусть  $V_4$ -- линейное пространство, состоящее из ростков формальных степенных рядов в начале координат вида
$$\Phi = \phi_4  +  \phi_5 +\dots = (f_3 +f_4+\dots, g_2 +g_3+\dots, h_4 +h_5+\dots)$$

В соответствии с теоремой 1 число параметров, от которых может зависеть отображение вида (\ref{nm}) $\Gamma$ на $\tilde{\Gamma}$ не превосходит размерности ${\rm Ker} \, L$ на пространстве $V_4$ .  Это, вместе с оценкой числа параметров в 3-струе, даст общую оценку числа прараметров, от которых может зависеть отображение и, в частности,  оценку размерности группы локальных автоморфизмов гиперповерхности $\Gamma$. Таким образом, для получения оценки размерности автоморфизмов 2-невырожденной гиперповерхности нам осталось дать оценку размерности ${\rm Ker} \, L$ на  $V_4$.

\vspace{5ex}

Оператор $L$ содержит большое число произвольных постоянных. Для того, чтобы упростить работу по оценке размерности ядра применим к уравнению $L(f,g,h)=0$ тот же самый прием, т.е. рекурсию на глубину два, но предварительно поменяв веса основных переменных. Зададим новые веса так:
$$[z]=[\bar{z}]=2, \;  [\zeta]=[\bar{\zeta}]=1, \;  [w]=[u]=4.$$
Если теперь, используя новое весовое разложение $\phi=(f,g,h)$, положить $\phi_{\mu}=(f_{\mu-2},g_{\mu-4},h_{\mu})$, то
${\mu}$-я весовая компонента $L(\phi)=0$ имеет вид
 \begin{eqnarray} \label{L3}
\nonumber   L_{\mu} =  2 \, {\rm Re} [ i \, h_{\mu} + i \,\Delta(h_{\mu-1})]+\qquad\qquad\qquad\qquad\qquad\qquad\\
\nonumber 2 \, {\rm Re} [ 2 \, <f_{\mu-2},\bar{z}> +  2 \,<\Delta (f_{\mu-3}),\bar{z}> +4 \, K(f_{\mu-3},z) \, \bar{\zeta}]+\qquad\qquad\qquad\qquad\\
\nonumber  2 \, {\rm Re} [ 2 \, \bar{K}( \bar{z}, \bar{z})\, g_{\mu-4} +2 \, \bar{K}( \bar{z}, \bar{z})\, \Delta(g_{\mu-5})+
(2 \, \bar{R}( \bar{z}, \bar{z}) \, \zeta + S(z,\bar{z}) \, \bar{\zeta}) \, g_{\mu-5}+\\(2\, \bar{P}( \bar{z}, \bar{z},z)+ \bar{Q}( \bar{z}, \bar{z},\bar{z}))\, g_{\mu-6}],  \quad   \mbox{ где   }     w=u+i \, <z,\bar{z}>=0.  \qquad\qquad\qquad
  \end{eqnarray}

\vspace{3ex}

{\bf Лемма 12:} Размерность пространства решений (\ref{L3}) не превосходит размерности пространства решений
$\mathcal{L}(f,g,h)=0$, где
 \begin{eqnarray} \label{L4}
\nonumber   \mathcal{L}(f,g,h) =  2 \, {\rm Re} [ i \, h + i \,\Delta(h)]+
2 \, {\rm Re} [ 2 \, <f,\bar{z}> +  2 \,<\Delta f,\bar{z}> +4 \, K(f,z) \, \bar{\zeta}]+\\
\nonumber  2 \, {\rm Re} [ 2 \, \bar{K}( \bar{z}, \bar{z})\, g +2 \, \bar{K}( \bar{z}, \bar{z})\, \Delta(g)+
2 \, \bar{R}( \bar{z}, \bar{z}) \, \zeta \, g + S(z,\bar{z}) \, \bar{\zeta} \, g],  \qquad\qquad\qquad \\
   \mbox{ где   }     w=u+i \, <z,\bar{z}>=0.  \qquad\qquad\qquad\qquad\qquad\qquad\qquad\qquad\qquad
  \end{eqnarray}
{\it Доказательство:}  Это сразу следует из теоремы 1.

\vspace{3ex}

Отметим при этом, что рекурсия, описанная оператором  $\mathcal{L}$, стартует с $\mu=5$. При этом нас интересует размерность ядра $L$ на пространстве $V_4$ в старой весовой градуировке. Поэтому лемма 12 нуждается в небольшой
коррекции.   Пусть $\tilde{V}_5$  состоит из наборов $(f,g,h)$, где $f=\tilde{O}(3), \; g=\tilde{O}(2), \; h=\tilde{O}(5)$  в соответствии с новым весом. Непосредственно убеждаемся в справедливости следующей леммы.

\vspace{3ex}

{\bf Лемма 13:}
   Если  $\phi=(f,g,h) \in V_4 \cap {\rm Ker} \,\mathcal{L}$, то $\phi \in \tilde{V}_5$. \\
{\it Доказательство:} Если $\phi \in V_4$, то $\phi = \chi +\psi$, где  $\psi \in \tilde{V}_5$, а $\chi=(0,0,\gamma \, \zeta^4)$.
Отделяя в соотношении  $\mathcal{L}(\chi +\psi)=0$ компоненту веса четыре, получаем  $\mathcal{L}(\chi )=0$. Откуда сразу следует, что $\chi=0$.

\vspace{3ex}

Переходя к оценке размерности  ядра оператора $\mathcal{L}$, отметим, что
оператор зависит от параметров $(k,m)$, ограничения на которые содержатся в лемме 5 (допустимые значения), и от трех коэффициентов квадратичной формы $R(z,z)=r_1 \, z_1^2 + r_2 \, z_1 \, z_2 + r_3 \, z_2^2$, которые не связаны никакими ограничениями. Также отметим, что независимо от значений параметров  ${\rm Ker} \,\mathcal{L}$ содержит двумерное
подпространство (тривиальные решения), которое, впрочем, не пересекается с  $\tilde{V}_5$.
\begin{eqnarray} \label{S1}
 (f_1=f_2=g=0, \quad h=t_1),   \;  t_1  \in \mathbf{R}, \\
\nonumber (f_1= t_2 \, z_1, \; f_2=t_2 \, z_2, \;  g=0, \quad h=t_2^2 \, w),   \;  t_2 \in \mathbf{R}.
\end{eqnarray}

\vspace{3ex}

{\bf Лемма 14:} Пусть $\phi=(f,g,h) \in \phi \in \tilde{V}_5 \cap {\rm Ker} \,\mathcal{L}$.\\
(a) Если $(k=1,m=0)$ (пара номер 9 из леммы 8) и $R(z,z)=r_1 \, z_1^2 $ , то\\
\begin{eqnarray*}
f_1=i \, \bar{n}_1 \, z_1^1, \quad  f_2=2 \, i \, \bar{n}_1 \, z_1 \, z_2 - \bar{n}_2 \, z_1^2 + n_1 \, w,\\
g=\frac{n_2 \, z_1-i \, n_1 \, z_2+ 2 \, i \, \bar{n}_1 \, z_1 \, \zeta}{1+2 \, \bar{r}_1 \, \zeta}, \quad h = 2 \, i \, \bar{n}_1 \, z_1 \, w,
\end{eqnarray*}
где $n_1$ и $n_2$ -- комплексные числа. Соответственно $\dim  ( \tilde{V}_5 \cap {\rm Ker} \,\mathcal{L}) = 4$.\\
(b) Во всех остальных случаях $\phi =0$.   Соответственно $\dim  ( \tilde{V}_5 \cap {\rm Ker} \,\mathcal{L}) = 0$.\\
{\it Доказательство:} Доказательство представляет собой рутинное, но объемное вычисление, которое осуществляется
средствами компьютерной алгебры (Maple).  Это вычисление проводится отдельно для пар с номерами
(1,2,3,4,5,6).  И отдельно для номеров (7,8,9). Для единообразного рассмотрения пар (1,2,3) и (4,5,6) вводится параметр
$\varepsilon=\pm 1$, учитывающий сигнатуру  формы Леви.
Введем также обозначения
\begin{eqnarray*}
(f_1(0,0,0,u),  f_2(0,0,0,u))=(a_1(u), a_2(u))=a(u), \\
g(0,0,0,u)=b(u), \;h(0,0,0,u)=c(u),\\
\frac{\partial \, f_1}{\partial \, z_1}(0,0,0,u)=a_{11}(u), \; \frac{\partial \, f_2}{\partial \, z_1}(0,0,0,u)=a_{21}(u), \\
\frac{\partial \, f_1}{\partial \, z_2}(0,0,0,u)=a_{12}(u), \; \frac{\partial \, f_2}{\partial \, z_2}(0,0,0,u)=a_{22}(u),\\
\frac{\partial \, g}{\partial \, z_1}(0,0,0,u)=b_{1}(u), \; \frac{\partial \, g}{\partial \, z_2}(0,0,0,u)=b_{2}(u),
 \; \frac{\partial^2 \, g}{\partial \, z_1^2}(0,0,0,u)=B(u).
\end{eqnarray*}
Схема вычисления в первом и втором случаях отличается мелкими деталями. Опишем её на примере второго случая -- пар (7,8,9).
{\it Шаг первый}. Положим в соотношении
\begin{equation}\label{M}
 \mathcal{L}(f_1,f_2,g,h)=0
\end{equation}
$\bar{z}=0, \; \bar{\zeta}=0,$  получим выражение $h(z_1,z_2,\zeta,u)$ через $(a_1(u),a_2(u),b(u),c(u))$. Это выражение имеет вид
\begin{eqnarray}  \label{H}
\nonumber h(z_1,z_2,\zeta,u)=\bar{c}(u)+2\,i\,<z,\bar{a}(u)>+ 2\,i\,\bar{b}(u) \,K(z,z) \qquad \qquad\qquad
\end{eqnarray}
При этом, подставляя $(z=0, \zeta=0)$, убеждаемся, что $\bar{c}(u)=c(u)$. \\
{\it Шаг второй}. Подставляем полученное значение $h$ в (\ref{M}), вычисляем $\mathcal{L}'_{\bar{z}_1}$ и
$ \mathcal{L}'_{\bar{z}_2}$, подставляем $\bar{z}=0, \; \bar{\zeta}=0$  и из этих соотношений получаем выражения для $f_1$ и $f_2$ через $(a_1(u),a_2(u),b(u),c(u),a_{11},a_{12},a_{21},a_{22})$. Они имеют вид
\begin{eqnarray*}
f_1=a_1(u)+2\,i\,\bar{b}'(u)\,m\,z_1\,z_2^2+2\,i\,\bar{b}'(u) \, z_1^3-4\,\bar{m}\,z_1\,\zeta\,\bar{b}(u)+c'(u)\,z_1+\\
2\,i\,\bar{a}'_2(u)\,z_1^2+2\,i\,\bar{a}'_1(u)\,z_1\,z_2-\bar{b}_2(u)\,m\,z_2^2-2\,\zeta\,\bar{a}_2(u)\,\bar{m}-\\
\bar{b}_2(u)\,z_1^2-\bar{a}_{22}(u)\,z_1-\bar{a}_{12}(u)\,z_2,\\
f_2=a_2(u)+2\,i\,\bar{b}'(u)\,m\,z_2^3+2\,i\,\bar{b}'(u)\,z_1^2\,z_2-\bar{b}_1(u)\,m\,z_2^2-4\,m\,z_2\,\zeta\,\bar{b}(u)+\\
c'(u)\,z_2+2\,i\,\bar{a}'_2(u)\,z_1\,z_2+2\,i\,\bar{a}'_1(u)\,z_2^2-\bar{b}_1(u)\,z_1^2-\\
\bar{a}_{11}(u)\,z_2-\bar{a}_{21}(u)\,z_1-2\,\zeta\,\bar{a}_1(u).
\end{eqnarray*}
 Подставляя $(z=0, \zeta=0)$, получаем
 $$a_{22}(u)=c'(u)-\bar{a}_{11}(u),  \quad {\rm Re} \, a_{12}(u)={\rm Re}\, a_{21}(u)=0.$$
 {\it Шаг третий.} Подставим полученные значения $f_1$ и $f_2$ в (\ref{M}), вычислим $\mathcal{L}''_{\bar{z}^2_1}$,  подставим $\bar{z}=0, \; \bar{\zeta}=0$ и из этого соотношения получим выражение для $g$ через $(a_1(u),a_2(u),b(u),c(u),a_{11},a_{12},a_{21},a_{22},b_1(u),b_2(u),B(u))$.
\begin{eqnarray*}
g=\frac{1}{2 \,(2 \,\bar{r}_1\,\zeta+1)}\;
(2\,b(u)+4\,i \,\bar{b}'_1(u)\,z_2\,z_1^2-\bar{B}(u)\,m\,z_2^2+12\,i\,\zeta\,\bar{a}'_1(u)\,z_2-2\,c'(u)\,\zeta+\\
4\,a_{22}(u)\,\zeta+
2\,b_1(u)\,z_1-\bar{B}(u)\,z_1^2+2\,b_2(u)\,z_2-8\,\bar{b}_1(u)\,m\,\zeta\,z_2-\\
2\,i\,\bar{a}'_{12}(u)\,z_2\,z_1+2\,i\,\bar{a}'_{12}(u)\,z_1\,z_2+
4\,i\,\bar{a}'_2(u)\,\zeta\,z_1+20\,i\,m\,z_2^2\,\zeta\,\bar{b}'(u)-2\,i\,{a}_{22}(u)\,z_2^2+\\
4\,i\,\bar{b}'(u)\,\zeta\,z_1^2+4\,\bar{b}''(u)\,z_2^4\,m+4\,\bar{b}''(u)\,z_2^2\,z_1^2+4\,\bar{a}'_2(u)\,z_2^2\,z_1+\\
2\,i\,\bar{a}'_{11}(u)\,z_2^2+4\,i\,\bar{b}'_1(u)\,z_2^3\,m+2\,\bar{a}''_1(u)\,z_2^3)
\end{eqnarray*}
  Подставляя $(z=0, \zeta=0)$, получаем $B=-\frac{1}{2}\,\bar{B}$, откуда следует, что $B=0$.\\
 {\it Шаг четвертый.} После подстановки в (\ref{M}) полученного выражения для $g$ мы получаем соотношение, которое
имеет вид вещественного полинома по $(z,\bar{z},\zeta,\bar{\zeta})$, коэффициенты которого суть дифференциальные полиномы от введенных функций переменного $u$ и их производных. Приравнивняем к нулю все коэффициенты. Анализ полученной системы обыкновенных дифференциальных уравнений позволяет завершить доказательство леммы.

\vspace{7ex}

{\bf 2-невырожденные гиперповерхности специального вида}

\vspace{3ex}

В соответствии с лемой 14  нетривиальное ядро имеется только для некоторого специального
класса 2-невырожденных гиперповерхностей. Т.е. таких, что в каждой своей точке
они могут быть заданы уравнением вида
$$v=2 \, {\rm Re} (z_1 \, \bar{z}_2) + 2 \, {\rm Re} (z_1^2 \, \bar{\zeta}) +2 \, {\rm Re} (r_1 \,z_1^2 \, \bar{\zeta}^2) +...$$
После замены $\zeta \rightarrow   \zeta+ \bar{r}_1 \, \zeta^2$ уравнение принимает вид
\begin{equation}\label{C}
v=2 \, {\rm Re} (z_1 \, \bar{z}_2) + 2 \, {\rm Re} (z_1^2 \, \bar{\zeta}) +\mbox{мономы нового веса 7 и выше}
\end{equation}

Для изучения таких гиперповерхностей нам будет удобно сделать перестановку координат и еще раз поменять веса. Пусть теперь
$$ [z_1]=2, \;  [z_2]=[\zeta]=1, \; [w]=[u]=3.$$
Тогда гиперповерхность (взвешенная модельная поверхность) задаётся соотношением
\begin{equation}\label{Q}
Q=\{v=2 \, {\rm Re} (z_1 \, \bar{\zeta} + z_2 \, \bar{\zeta}^2) \}
\end{equation}
при таком выборе весов  $Q$ -- это график квазиоднородного вещественного полинома веса 3.   Это позволяет
применить рекурсию на глубину один и получить исчерпывающий ответ.   Отметим также, что использование взвешенных модельных поверхностей применяется достаточно давно (\cite{VB96}, \cite{AE01}, \cite{KMZ14}) и эта техника вполне стандартна.

Подгруппа $\mathcal{Q}$ автоморфизмов гиперповерхности $Q$, которая обеспечивает голоморфную однородность $Q$, состоит из преобразований вида
\begin{eqnarray}\label{QQ}
\nonumber
z_1 \rightarrow  a + z_1 ,  \qquad z_2 \rightarrow b + 2 \, \bar{a} \, \zeta + z_2, \qquad \zeta  \rightarrow  c + \zeta , \qquad\qquad \\
w    \rightarrow d  + 2 \, i \, (a \, \bar{b}  + a^2 \, \bar{c} +  (\bar{b} +2 \, a \, \bar{c})\, z_1 + \bar{a} \, z_2   +
\bar{a}^2 \, \zeta + \bar{c} \, z_1^2) + w,
\end{eqnarray}
где $(a,b,c,d)$ -- произвольная точка $Q$.

Пусть   $\Gamma_{0}$  -- росток гиперповерхности в начале координат вида
\begin{equation}\label{GS}
v=2 \, {\rm Re} (z_1 \, \bar{\zeta} + z_2 \, \bar{\zeta}^2) +O(4),
\end{equation}
где $O(4)$ - это слагаемые веса четыре и выше.
Рассмотрим отображение $\phi=(f, g,h,e)$ этого ростка на другой росток такого же вида.
Причем
\begin{equation}\label{MG}
f=z_1+f_3+\dots, \; g=z_2+g_2+\dots, \; h=\zeta+h_2+\dots, \;e=w+e_4+\dots
\end{equation}
  (нижние индексы обозначают веса компонент).  Тогда, записывая в виде аналитического соотношения тот факт, что это отображение переводит первую гипеповерхность во вторую и отделяя в нем $\mu$-ю весовую компоненту, получаем
\begin{eqnarray}
\nonumber
- \,{\rm Im}\, e_{\mu} + 2 \,{\rm Re}( f_{\mu-1}\, \bar{\zeta} +g_{\mu-2}\, \bar{\zeta}^2 +h_{\mu-2} \,(\bar{z}_1 +
2 \, \bar{z}_2 \,\zeta)) =\dots
\end{eqnarray}
где $w=u+2 \, i  \,  {\rm Re} (z_1 \, \bar{\zeta} + z_2 \, \bar{\zeta}^2) \}$, а многоточие означает выражение, зависящее от компонент с меньших весов (т.е. для $f$ -- меньше $\mu-1$, для $g$ и $h$ -- меньше $\mu-2$, для $e$ -- меньше $\mu$.

Таким образом мы видим, что размерность семейства отображений вида (\ref{MG}) контролируется размерностью ядра гомологического оператора.
\begin{eqnarray} \label{OP}
 L(f,g,h,e)= 2 \,{\rm Re} (i \, h + 2 \, f \, \bar{\zeta} + 2 \, g \, \bar{\zeta}^2 + 2 \,h \,(\bar{z}_1 +2 \, \bar{z}_2 \,\zeta))
 \end{eqnarray}
при  $w=u+2 \, i  \,  {\rm Re} (z_1 \, \bar{\zeta} + z_2 \, \bar{\zeta}^2) $.

С другой стороны, если
\begin{eqnarray*}
X=2 \, {\rm Re} \left( f  \,\frac{\partial}{\partial z_1}+g \,\frac{\partial}{\partial z_2}+
h \,\frac{\partial}{\partial \zeta}+e \,\frac{\partial}{\partial w} \right)
\end{eqnarray*}
-- росток векторного поля в начале координат, т.ч. $(f, g,h,e)$ - голоморфно в нуле, тогда
равенство $L(f, g,h,e)=0$ равносильно тому, что $X $ -- элемент алгебры Ли инфинитезимальных автоморфизмов
$Q$ в начале координат -- ${\rm aut} \, Q$.

Веса, введенные нами для координат пространства, естественно продолжаются и на дифференцирования по этим
координатам. Дифференцирование по $z_1$ имеет вес (-2), по $z_2$ и по $\zeta$  -- вес (-1), по $w$ -- вес (-3).
Это превращает ${\rm aut} \, Q$ в градуированную алгебру Ли вида $g_{-3}+g_{-2}+\dots$. Подалгебра $g_0$
 содержит градуирующее поле
 $$X_0 =2 \, {\rm Re} \left( 2 \, \frac{\partial}{\partial z_1}+ 1 \, \frac{\partial}{\partial z_2}+
 1 \, \frac{\partial}{\partial \zeta}+3 \,\frac{\partial}{\partial w}\right).$$
В такой ситуации если некоторое поле есть элемент алгебры, то каждая его градуированная компонента -- тоже.
Рассуждение из работы В.Каупа \cite{K73},  позволяет утверждать, что  алгебра ${\rm aut} \, Q$,
в таком случае, обязана быть конечно градуированной (полиномиальной).  Но мы не будем использовать это утверждение, а вычислим алгебру явно.

\vspace{3ex}

Переходим к  вычислению алгебры ${\rm aut} \, Q$, которая совпадает с ядром оператора (\ref{OP}).  Процедура
вычисления аналогична той, что описана в доказательстве леммы 14.  Однако само вычисление  проще.

\vspace{3ex}

Введем обозначения
\begin{eqnarray*}
f(0,0,0,u) =a(u), \;  g(0,0,0,u)=b(u), \;h(0,0,0,u)=c(u), \qquad \qquad\\
e(0,0,0,u)=d(u),  \; \frac{\partial \, f}{\partial \, z_1}(0,0,0,u)=a_{1}(u), \frac{\partial \, f}{\partial \, \zeta}(0,0,0,u)=a_{3}(u), \qquad\\
\frac{\partial \, g}{\partial \, z_1}(0,0,0,u)=b_{1}(u), \; \frac{\partial \, g}{\partial \, \zeta}(0,0,0,u)=b_{3}(u), \qquad\qquad\qquad\qquad\qquad\qquad\\
\frac{\partial \, h}{\partial \, z_1}(0,0,0,u)=c_{1}(u), \; \frac{\partial \, h}{\partial \, \zeta}(0,0,0,u)=c_{3}(u),
\frac{\partial^2 \, g}{\partial \, z_1^2}(0,0,0,u)=B(u).
\end{eqnarray*}
 Положим в соотношении
\begin{equation}\label{QL}
 \mathcal{L}(f,g,h,e)=0
\end{equation}
$\bar{z}_1=0, \; \bar{z}_2=0, \; \bar{\zeta}=0,$  получим  выражение для $h$. Это полином степени два от
$(z_1,z_2,\zeta)$ с коэффициентами,  зависящими от  $(a(u),b(u),c(u),d(u))$. Подставляя это значение $h$ в (\ref{QL}), вычисляем $\mathcal{L}'_{\bar{\zeta}}$, подставляем $\bar{z}=0, \; \bar{\zeta}=0$  и из найденного соотношения получаем  выражение для $f$, которое является полиномом степени  три с коэффициентами,  зависящими от  $(a,a',b',c',d,a_3,b_3,c_3)$. Вычисляем $\mathcal{L}'_{\bar{z}_1}$, подставляем $\bar{z}=0, \; \bar{\zeta}=0$  и из этого соотношения получаем  выражение для $h$, которое является полиномом степени
два  с коэффициентами,  зависящими от $(a,a',b',c,c',a_1,b_1,c_1)$ . Подставляя эти значения $f$ и $h$ в (\ref{QL}), вычисляем $\mathcal{L}''_{\bar{\zeta}^2}$, положим $\bar{z}=0, \; \bar{\zeta}=0$  и из найденного соотношения получим  выражение для $g$, которое является полиномом степени четыре с коэффициентами,  зависящими от  $(a',a'',b,b',b'',c',c'',d',a'_1,b_1,b_3,b'_1,c'_1,c'_3,B)$.

Дальнейший анализ соотношения (\ref{QL}) дает
\begin{eqnarray*}
{\rm Im} \, d={\rm Re} \, c_1={\rm Re} \, B=0, \; b_3=i \, c',  \; c_3=d'-\bar{a}_1, \\
 a'=b'=c''=d''=a'_1=a_3=b'_1=c'_1=B'=0.
\end{eqnarray*}

Подсчитаем число сободных вещественных параметров
$$a   -   2, \; b  -  2, \; c  -  4, \; d  -  2, a_1  -  2, \; a_3  -  0, \, b_1  -  2, \; b_3  -  0, \; c_1  -  1,  \; c_3   -   0, \; B  -  1.$$
Таким образом, получаем, что размерность ${\rm aut} \, Q$ не превосходит 16.

\vspace{3ex}
С другой стороны не трудно выписать несколько младших весовых компонент ${\rm aut} \, Q$. Вот эти компоненты:
\begin{eqnarray} \label{DQ}
  g_{-3}=\{(0,\;0,\;0,\;d) \},    \quad  d \in \mathbf{R}, \qquad \qquad \qquad \qquad \qquad \qquad \\
\nonumber
  g_{-2}=\{(a,\;0,\;0,\;2 \,i \, \bar{a}\, \zeta) \},    \quad  a \in \mathbf{C}, \qquad \qquad \qquad \qquad \qquad\\
\nonumber
g_{-1}=\{(- 2 \, \bar{c} \, z_2 + i \, k \, \zeta ,\; b, \; c, \; 2 \,i \, \bar{c}\, z_1+2 \, i \, \bar{b} \, \zeta^2) \},    \quad  b, \; c \in \mathbf{C}, \;  k \in \mathbf{R}, \qquad \qquad\qquad\qquad\\
\nonumber
g_{0}= \{ (M \,  z_1 + m \, \zeta^2 \, ( 2 \, M - 3 \, l) \, z_2 - \bar{m} \, \zeta \,   (3 \, l - \bar{M})\, \zeta, 3 \,l \, w)  \},    \;
M, \; m \in \mathbf{C}, \;  l \in \mathbf{R}, \\
\nonumber
g_{1}=\{(2\,i\,\bar{N}\,z_1\,\zeta+N\,w,   2 \,i\, \bar{N}\, z_2\,\zeta -i\,N\, z_1+i\,n \, \zeta^2 , i\,\bar{N}\,\zeta^2,  2\,i\,\bar{N}\,\zeta w) \},    \quad  N  \in \mathbf{C}, \; n \in \mathbf{R}.
\end{eqnarray}
И мы видим, что размерность суммы этих пяти компонент равна 16.   Таким образом алгебра вычислена.
Для дальнейшего нам будет удобно представить $g_{-1}$  в виде прямой суммы $g'_{-1}+g''_{-1}$, где
\begin{eqnarray*}
g'_{-1}=\{(- 2 \, \bar{c} \, z_2  ,\; b, \; c, \; 2 \,i \, \bar{c}\, z_1+2 \, i \, \bar{b} \, \zeta^2) \},\\
g''_{-1}=\{(  i \, k \, \zeta,\;0,\;0,\;0) \}.
\end{eqnarray*}

Сформулируем полученный результат.

{\bf Теорема 15:}\\
(a) Алгебра ${\rm aut} \, Q$ -- это сумма пяти градуированных компонент $g_{-3}+g_{-2}+g_{-1}+g_{0}+g_{1}$, сами
компоненты выписаны выше (\ref{DQ}).   $\dim \, {\rm aut} \, Q =16$. \\
(b) При этом ${\rm aut}_0 \, Q$, стабилизатор начала координат в ${\rm aut} \, Q$ (т.е. поля из алгебры, обращающиеся в ноль в начале координат) -- это $g''_{-1}+g_0+g_1$, его размерность равна 9.\\
 (c) Подалгебра $g_{-3}+g_{-2}+g'_{-1}$ -- это алгебра Ли подгруппы $\mathcal{Q}$ -- группы $"$сдвигов$"$.  Причем $Q$ находится в естественном взаимно-однозначном соответствии с $\mathcal{Q}$.  Это позволяет  перенести на  $\mathcal{Q}$ структуру  вложенной гиперповерхности.\\
 (d)  Если $\Gamma_0$ - росток гиперповерхности вида (\ref{DQ}), имеет место оценка как для всей алгебры, так и отдельно для стабилизатора начала координат .
$$ \dim \, {\rm aut} \,\Gamma_0 \leq 16, \quad   \dim \, {\rm aut}_0 \,\Gamma_0 \leq 9.$$

\vspace{3ex}

Для полноты картины выпишем автоморфизмы, порожденные этими полями. Как было отмечено, поля из
$g''_{-1}+g_0+g_1$ порождают $\mathcal{Q}$.

Поле $(i \, \zeta,0,\;0,\;0)$  из $g'_{-1}$ порождает преобразование
$$z_1 \rightarrow z_1 +i \, t \, \zeta, \; z_2 \rightarrow  z_2,   \zeta \rightarrow \zeta , \; w  \rightarrow w.$$
Аналогично поле   $(0, \; i \, \zeta^2,\;0,\;0)$  из $g_{1}$ порождает преобразование
$$z_1 \rightarrow  z_1, \;z_2 \rightarrow z_2 +i \, t \, \zeta^2, \;   \zeta \rightarrow \zeta , \; w  \rightarrow w.$$

Для вычисления преобразований, порожденных $g_0$, положим $M=2 \, -\mu/3$. Тогда получаем
\begin{eqnarray*}
z_1  \rightarrow  \left( z_1 + m \, \left(\frac{e^{\mu t}-1}{\mu}\right) \, \zeta^2  \right) \, e^{(2l-\mu/3)t}, \\
z_2  \rightarrow  \left( z_2 - \bar{m} \, \left(\frac{e^{\mu t}-1}{\mu}\right) \, \zeta  \right) \, e^{(l- 2\mu/3)t},\\
\zeta \rightarrow \zeta \, e^{(l+ \mu /3)t}, \quad w \rightarrow w \, e^{3lt}. \qquad\qquad
\end{eqnarray*}
И, наконец преобразования из $g_1$ при $n=0$ имеют вид
\begin{eqnarray*}
z_1 \rightarrow\frac{z_1}{(1-i\, \bar{N} \, \zeta \, t)^2},   \quad z_2   \rightarrow\frac{z_2-i\,N\, z_1\,t}{(1-i\, \bar{N} \, \zeta \, t)^2},\\
\zeta \rightarrow \frac{\zeta}{1-i\, \bar{N} \, \zeta \, t}, \quad    w \rightarrow \frac{w}{(1-i\, \bar{N} \, \zeta \, t)^2}.
\end{eqnarray*}

\vspace{3ex}
Гипеповерхность $Q$ замечательна во многих отношениях.  Она представляет собой общее начало двух последовательностей гипеповерхностей пространства $\mathbf{C}^N$ при $N \geq 4$ (больше у них персечений нет) .
Первая последовательность была рассмотрена в работе А.Лабовского (\cite{AL97}, 1997) как пример голоморфно однородных $l$-невырожденных гиперповерхностей с произвольным $l$. Если эта гиперповерхность расположена
в $\mathbf{C}^N$, то она равномерно $(N-2)$-невырождена.

С другой стороны,  в недавней работе И.Зеленко и Д.Сайкса \cite{ZS21}  была описана серия голоморфно однородных 2-невырожденных гиперповерхностей пространства $\mathbf{C}^N$ с алгеброй голоморфных автоморфизмов размерности $(N-1)^2+7$ и доказано, что эти гиперповерхности оптимальны в классе голоморфно однородных.
Т.е. никакая голоморфно однородная гиперповерхность не может иметь автоморфизмы большей размерности.
Отметим, что эта работа использует технику, весьма далекую от нашей. Это дифференциальная геометрия в стиле Э.Картана и Н.Танаки.

\vspace{3ex}

{\bf Утверждение 16:} (a) Если вещественная гиперповерхность $\Gamma$ всюду, кроме собственного аналитического подмножества, является 2-невырожденной, то для любой ее точки $\xi$ размерность алгебры автоморфизмов ростка гиперповерхности в этой точке ${\rm aut} \,\Gamma_{\xi}$ не превосходит 17.\\
(b)  Если же эта гиперповерхность в точке общего положения   принадлежит специальному классу (\ref{C})
(или более широкому классу  (\ref{GS})), то $\dim \, {\rm aut} \,\Gamma_{\xi} \leq 16$.\\
{\it Доказательство:} Размерность   ${\rm aut} \,\Gamma$ не превосходит размерности гиперповерхности, которая
равна 7 плюс размерность стабилизатора точки. Для оценки стабилизатора в произвольной точке достаточно провести оценку в точке 2-невырожденности. В соответствии с теоремой 1 и всеми последующими леммами размерность стабилизатора оценивается через размерность младшей струи (леммы 4 и 8) и размерность ядра $\mathcal{L}$  на $\tilde{V}_5$ (которая равна нулю).   Размерность группы параметров $(C,\rho,\beta)$ не превышает 3. Прараметры $(a, \alpha)$ дают еще 8. Итого $7+ 3+8=18$. Однако, что бы получить 18 надо, чтобы орбита начала координат была 7-мерной. Это означает голоморфную однородность. Но тогда, в соответствии с \cite{ZS21},  размерность не выше 16-ти.
Поэтому мы можем считать, что размерность орбиты меньше 7. Откуда получаем оценку $6+3+8=17$.\\
В случае (b) мы можем воспользоваться теоремой 15 (d).  Утверждение доказано.

\vspace{3ex}

{\bf Теорема 17:} Пусть $\Gamma$ -- голоморфно невырожденная вещественно аналитическая гиперповерхность в $\mathbf{C}^4$,  точка $\xi \in \Gamma$ и $\Gamma_{\xi}$ -- росток  $\Gamma$ в точке $\xi$. Пусть ${\rm aut} \,\Gamma_{\xi}$ -- алгебра Ли инфинитезимальных голоморфных автоморфизмов ростка. Тогда:\\
(1)
$$  \dim \, {\rm aut} \,\Gamma_{\xi} \leq 24.$$
(2)  Если известно, что $\Gamma$  2-невырождена всюду, кроме собственного аналитического подмножества, то можно утверждать, что
$$\dim \, {\rm aut} \,\Gamma_{\xi} \leq 17.$$
(3)  Если известно, что $\Gamma$  3-невырождена всюду, кроме собственного аналитического подмножества, то можно утверждать, что
$$\dim \, {\rm aut} \,\Gamma_{\xi} \leq 20.$$

\vspace{3ex}

{\bf Теорема 18:} Пусть $\Gamma_{\xi}$ -- росток произвольной вещественно аналитической гиперповерхности в $\mathbf{C}^4$ т.ч  $\dim \, {\rm aut} \,\Gamma_{\xi} = 24$,   тогда $\Gamma_{\xi}$ эквивалентен одной из двух стандартных невырожденных гиперквадрик (т.е. Леви-невырожден и сферичен).\\
{\it Доказательство:} Если $\Gamma_{\xi}$ не является Леви-невырожденной в общей точке, то, как следует из теоремы 17, размерность не превосходит 20. Таким образом, $\Gamma_{\xi}$  Леви-невырождена в общей точке. Если она там не сферична, то,  как доказано  в  \cite{KR15},  размерность не превосходит 13-ти. Поэтому она сферична.
Но тогда, в соответствии с другим результатом Б.Кругликова \cite{KR20}, если  $\Gamma_{\xi}$ не эквивалентна
гиперквадрике (произвольной сигнатуры), то размерность не выше 17. Теорема доказана.

\vspace{3ex}

Аналогичные оценки для $\mathbf{C}^2$ и $\mathbf{C}^3$  -- это 8 и 15. Они также достигаются  только на гиперквадриках. Проблемными, как и в  $\mathbf{C}^4$, являются гиперповерхности, которые в общей точке являются сферическими.   Результат для $\mathbf{C}^2$ -- это работа  И.Коссовского и Р.Шафикова \cite{KS}, а для $\mathbf{C}^3$ -- А.Исаева и Б.Кругликова \cite{IK}.

\vspace{3ex}

Эти результаты вместе с известным критерием конечномерности дают следующий список возможностей. Пусть
$d=\dim \, {\rm aut} \,\Gamma_{\xi}$. Тогда\\
\\
(1) $d=\infty$ тогда и только тогда, когда $\Gamma$  голоморфно вырождена.\\
(2) $d=24$ тогда и только тогда, когда $\Gamma$  эквивалентна одной из двух невырожденных стандартных гиперквадрик.\\
(3) Если $\Gamma_{\xi}$  несферичен, но $\Gamma$  сферична в точке общего положения, то $d \leq 17$.\\
(4)  Если $\Gamma$  в точке общего положения 1-невырождена (Леви-невырождена) и несферична, то $d \leq 13$.\\
(5)  Если $\Gamma$ -- 2-невырождена в точке общего положения, то $d \leq 17$.\\
(6)  Если $\Gamma$  -- 2-невырождена в точке общего положения и однородна (в окрестности $\xi$), то $d \leq 16$.\\
(7)  Если $\Gamma$  -- 3-невырождена в точке общего положения, то $d \leq 20$.\\

\vspace{3ex}

В этом списке пункты (1) и (2) вполне конкретны. То же следует сказать и о пункте (3). Действительно,  в работе \cite{KR20} имеется пример гиперповерхности такого типа, для которой оценка 17 реализуется. То же самое касается и пункта (4).  В \cite{KR15} имеется пример такой гиперповерхности, у которой размерность алгебры автоморфизмов равна 13. Пункт (6) также точен. Как было показано, он подкреплен примером. Поэтому для точной оценки из пункта (5) есть ровно две возможности -- 16 и 17. Последний пункт (7) самый неопределенный. Пока можно сказать только, что максимум не меньше 8 и не больше 20.

\vspace{3ex}

{\bf Вопрос  19:}  Каковы точные значения  максимумов  из пунктов  (5) и (7)?

\vspace{3ex}

{\bf Вопрос  20:}   Верно ли, что альтернатива остается верной для гиперповерхностей в размерности выше, чем 4?
Т.е. либо бесконечность, либо не больше, чем у гиперквадрики, у которой в $\mathbf{C}^N$ размерность  $(N+1)^2-1$.

У этого довольно старого  вопроса \cite{VB96} есть более свежая  версия (\cite{VB20}, conjecture (5.a),(5.b)).


\begin{thebibliography}{30}
\bibitem{P}   H.Poincare, Les fonctions analytiques de deux variables et la representation conforme,
Rend. Circ. Mat., Palermo, 1907, pp.185 -- 220.
\bibitem{VB05}  В.К.Белошапка, “Симметрии вещественных гиперповерхностей трехмерного комплексного пространства”, Матем. заметки, т.78, № 2 (2005),  сс.171–179.
\bibitem{BER} M.S.Baouendi, P.Ebenfelt, L.P.Rothschild, CR automorphisms of real analytic manifolds in complex space, Comm. Anal. Geom. 6 (1998), no. 2, 291–315.
\bibitem{KS} I.Kossovskiy, R.Shafikov, Analytic differential equations and spherical real hypersurfaces.
J. Differential Geom. 102 (2016), no. 1, pp.67 --126.
\bibitem{IK}  A.Isaev, B. Kruglikov,  On the symmetry algebras of 5-dimensional CR-manifolds. Adv. Math. 322 (2017), 530--564.
\bibitem{VB96} V.K.Beloshapka, Automorphisms of Degenerate Hypersurfaces in $\mathbf{C}^2$ and a Dimension Conjecture // Russian Journal of Mathematical Physic, 1996, vol.4, no.3, p.393 -- 396.
\bibitem{VB20} V.K.Beloshapka,    CR-Manifolds of Finite Bloom–Graham Type: the Method of Model Surface // Russian Journal of Mathematical Physics vol. 27, no.2, pp.155–174 (2020).
\bibitem{AL97}  А.С.Лабовский, О размерности группы биголоморфных автоморфизмов вещественно-аналитических гиперповерхностей, Матем. заметки, 61:3 (1997),  349–358.
 \bibitem{AE01}   А.Е.Ершова, Автоморфизмы 2-невырожденных гиперповерхностей в $\mathbf{C}^3$,
 Матем. заметки, 69:2 (2001),  214–222.
 \bibitem{KMZ14}    Kolar, M., Meylan, F., Zaitsev, D., Chern-Moser operators and polynomial models in CR geometry, Adv. Math. 263 (2014), 321–356.
 \bibitem{K73}  W. Kaup, Einige Bemerkungen uber polynomiale Vektorfelder,
Jordanalgebren und die Automorphismen von Siegelschen Gebieten,  Math. Ann. 204, 131 --144 (1973).
\bibitem{ZS21} D. Sykes and I. Zelenko,  Maximal Dimension of Groups of Symmetries of
Homogeneous 2-nondegenerate CR-structures of Hypersurface Type with a 1-dimensional Levi Kernel, 2021,  https://arxiv.org/abs/2102.08599.
 \bibitem{KR20} B. Kruglikov, Blow-ups and infinitesimal automorphisms of CR-manifolds.
Math. Z. 296 (2020), no. 3-4, 1701–1724.
 \bibitem{KR15} B. Kruglikov, Submaximally symmetric CR-structures. J. Geom. Anal. 26, 3090–3097 (2016).
 \bibitem{AS17} A.Santi,  Homogeneous models for Levi-degenerate CR manifolds, Kyoto J.Math. 60(1), 291–334 (2020).
\bibitem{FK08} G.Fels, W.Kaup,   Classification of Levi degenerate homogeneous CR-manifolds in dimension 5,
Acta Math., 201 (2008), 1–82.
\end{thebibliography}
\end{document}